\documentclass[a4paper]{amsart}

\usepackage{amsmath,amsthm,amsfonts,amssymb,amscd}
\input{xy}
\xyoption{all}
\newtheorem{theo}{Theorem}[section]
\newtheorem{pro}[theo]{Proposition}
\newtheorem{lem}[theo]{Lemma}

\newtheorem{defi}[theo]{Definition}
\newtheorem{rem}[theo]{Remark}

\newenvironment{sis}{\left\{\begin{aligned}}{\end{aligned}\right.}

\numberwithin{equation}{section}

\newcommand{\N}{\mathbb{N}}

\newcommand{\0}{\underline{0}}

\newcommand{\g}{\mathfrak{g}}

\newcommand{\ep}{\epsilon}

\renewcommand{\d}{{\rm d}}

\newcommand{\Sq}{{\rm Sq}}

\begin{document}

\title{Infinitesimal deformations of restricted simple Lie algebras II}
\author{Filippo Viviani}

\date{02 May 2008}

\address{Institut f\"ur Mathematik, Humboldt Universit\"at zu Berlin, 10099 Berlin (Germany).}
\email{{\tt viviani@math.hu-berlin.de}}



\thanks{During the preparation of this manuscript, 
the author was supported by a grant from the Mittag-Leffler Institute
(Stockholm).}

\maketitle

\begin{abstract}
We compute the infinitesimal deformations of two families of restricted simple
modular Lie algebras of Cartan-type: the Contact
and the Hamiltonian Lie algebras.
\end{abstract}

\section{Introduction}

Simple Lie algebras over an algebraically closed field of positive characteristic 
different from $2$ and $3$ were classified by Wilson-Block (see \cite{BW}) in the 
restricted case and by Strade (see \cite{STR6}) and Premet-Strade (see \cite{PS3}) in
the general case.
The classification remains still open in characteristic $2$ and $3$
(see \cite[page 209]{STR}).

According to this classification, \emph{simple modular} (that is over a field of positive
characteristic) Lie algebras are divided into two big families, called classical-type
and Cartan-type algebras.
The algebras of classical-type are obtained by the simple Lie algebras
in characteristic zero by first taking a model over the integers (via Chevalley bases)
and then reducing modulo $p$ (see \cite{SEL}). The algebras of Cartan-type were constructed by
Kostrikin-Shafarevich in 1966 (see \cite{KS}) as finite-dimensional analogues
of the infinite-dimensional complex simple Lie algebras, which occurred in Cartan's classification
of Lie pseudogroups, and are divided into four families, called Witt-Jacobson, Special,
Hamiltonian and Contact algebras. The Witt-Jacobson Lie algebras
are derivation algebras of truncated divided power algebras and the remaining three families
are the subalgebras of derivations fixing a volume form, a Hamiltonian form and a contact form,
respectively. Moreover in characteristic $5$ there is one exceptional simple modular Lie algebra
called the Melikian algebra (introduced in \cite{MEL}).

A particular important class of simple modular Lie algebras are the ones
which are \emph{restricted}. These can be characterized as those
modular Lie algebras such that the $p$-power of an inner derivation
(which in characteristic $p$ is a derivation) is still inner (see \cite{FS}
or \cite{STR}).
Important examples of restricted Lie algebras are the ones coming
from groups schemes. Indeed, there is a 
bijection between restricted Lie algebras over $k$ and infinitesimal
$k$-group schemes of height one (see \cite[Chap. 2]{DG}).

This paper is devoted to the study of the \emph{infinitesimal deformations}
of the restricted simple Lie algebras.  
The simple Lie algebras of classical type are known to be rigid
over a field of characteristic different from $2$ and $3$ (see \cite{RUD}),
in analogy of what happens in characteristic zero. In the papers \cite{Viv} and \cite{Viv2}, 
the author computed the infinitesimal deformations of the Witt-Jacobson, Special and Melikian 
restricted simple Lie algebras. In this paper, we compute the infinitesimal deformations
of the Contact algebras $K(n)$ and the Hamiltonian algebras $H(n)$ over a field
$F$ of characteristic different from $2$ and $3$.

By standard facts of deformation theory, the infinitesimal deformations of a
Lie algebra are parametrized by
the second cohomology of the Lie algebra with values in the adjoint representation
(see for example \cite{GER1}).

Before stating the main results of this paper, we recall that 
there is a canonical way to produce $2$-cocycles in $Z^2(\g,\g)$ for a modular Lie algebra $\g$ 
over a field of characteristic $p>0$, namely the squaring operation (see \cite{GER1}).
Given an element $\gamma\in \g$, one defines the squaring of $\gamma$ to be
\begin{equation}\label{Square}
\Sq(\gamma)(x,y)=\sum_{i=1}^{p-1}\frac{[{\rm ad}(\gamma)^i(x),{\rm ad}(\gamma)^{p-i}(y)]}
{i!(p-i)!}\in Z^2(\g,\g)
\end{equation}
where ${\rm ad}(\gamma)^i$ is the $i$-iteration of the inner derivation ${\rm ad}(\gamma)$. 

Assuming the (standard) notations from sections 2.1 and 3.1
about the Contact algebras $K(n)$ and the Hamiltonian algebras $H(n)$, 
we can state the main results of this paper.

\begin{theo}\label{K-finaltheorem}
Let $n=2m+1\geq 3$. Then
$$H^2(K(n),K(n))=\bigoplus_{i=1}^{2m} \langle {\rm Sq}(x_i) \rangle_F \oplus 
\langle {\rm Sq}(1)\rangle_F.$$
\end{theo}

 \begin{theo}\label{H-finaltheorem}
Let $n=2m\geq 2$. Then if $n\geq 4$ we have that
$$H^2(H(n),H(n))=\bigoplus_{i=1}^n \langle {\rm Sq}(x_i)\rangle_F 
\bigoplus_{\stackrel{i<j}{j\neq i'}} \langle \Pi_{ij} \rangle_F \bigoplus_{i=1}^m 
\langle \Pi_i \rangle_F \bigoplus \langle \Phi\rangle_F,$$
where the above cocycles are defined (and vanish outside) by 
$$\begin{sis}
& \Pi_{ij}(x^a, x^b)=x_{i'}^{p-1}x_{j'}^{p-1}[D_i(x^a)D_j(x^b)-D_i(x^b)D_j(x^a)] 
\hspace{0,3cm} \text{ for } j\neq i,i',  \\
& \Pi_i(x_ix^a,x_{i'}x^b)=x^{a+b+(p-1)\ep_i+(p-1)\ep_{i'}} \hspace{2,7cm} 
\text{ if } a+b<\sigma^i,\\
& \Pi_i(x_k,x^{\sigma^i})=-\sigma(k)x^{\sigma-\ep_{k'}} \hspace{3,3cm} 
\text{ for any }\hspace{0,3cm} 1\leq k\leq n,\\
& \Phi(x^a,x^b)=\sum_{\stackrel{\0<\delta\leq a, \widehat{b}}{|\delta|=3}}
\binom{a}{\delta}\binom{b}{\widehat{\delta}}\sigma(\delta)\:\delta!\:
x^{a+\widehat{b}-\delta -\widehat{\delta}}.\\
\end{sis}$$
If $n=2$ then 
$$H^2(H(2),H(2))=\bigoplus_{i=1}^2 \langle {\rm Sq}(x_i)\rangle_F 
\bigoplus \langle \Phi\rangle_F.$$

\end{theo}

In a forthcoming paper \cite{Viv3}, we use the above computations to determine the infinitesimal 
deformations of the simple finite group schemes in positive characteristic 
associated to the restricted simple Lie algebras.

The result presented here constitute part of my doctoral thesis. I thank my advisor
prof. Schoof for useful advices and constant encouragement. I thank the referee for
the suggestion of using the results of \cite{Far} and \cite{FS91} in order to simply 
some of the proofs of this paper.

\section{Contact algebra}

\subsection{Definition and Basic properties}

We first introduce some notations about the set $\N^n$ of $n$-tuple of natural numbers. 
We consider the order relation defined by
$a=(a_1,\cdots, a_n)<b=(b_1,\cdots,b_n)$ if $a_i<b_i$ for every $i=1,\cdots, n$.  
We define the degree of $a\in \N^n$ as $|a|=\sum_{i=1}^n a_i$ and the 
factorial as $a!=\prod_{i=1}^n a_i!$. For two multindex $a, b\in \N^n$ such that $b\leq a$, we 
set 
$\binom{a}{b}:=\prod_{i=1}^n \binom{a_i}{b_i}=\frac{a!}{b!(a-b)!}.$
For every integer
$j\in \{1,\cdots,n\}$ we call $\ep_j$ the $n$-tuple having $1$ at the $j$-th entry 
and $0$ outside.

Throughout this section we fix a field $F$ of characteristic $p\neq 2,3$ and 
an odd integer $n=2m+1\geq 3$. 
For any $j\in\{1,\cdots,2m\}$, we define the sign 
$\sigma(j)$ and the conjugate $j'$ of $j$ as follows:
$$\sigma(j)=\begin{sis}
1& \quad \text{ if }1\leq j\leq m,\\
-1& \quad \text{ if } m< j\leq 2m,\\
\end{sis}
\hspace{0,5cm} \text{ and } \hspace{0,5cm}
j'=\begin{sis}
j+m& \quad \text{ if } 1\leq j\leq m,\\
j-m& \quad \text{ if } m<j\leq 2m.
\end{sis}$$
Given a multindex $a=(a_1,\cdots,a_{2m})\in \N^{2m}$, we define the sign of $a$ as
$\sigma(a)=\prod\sigma(i)^{a_i}$ and the conjugate of $a$ as the multindex 
$\hat{a}$ such that $\hat{a}_i=a_{i'}$ for every $1\leq i\leq 2m$.
We are going to use often the following special $n$-tuples: 
$\0:=(0,\cdots, 0)$, $\tau:=(p-1,\cdots,p-1)$ and $\sigma:=(p-1,\cdots,p-1,0)$.

Let $A(n)=F[x_1,\cdots,x_n]/(x_1^p,\cdots, x_n^p)$ the ring of $p$-truncated 
polynomials in $n$-variables. Note that $A(n)$ is a finite $F$-algebra 
of dimension $p^n$ with a basis given by the monomials $\{x^a=x_1^{a_1}\cdots x_n^{a_n} 
\:|\: a\in \N^n, \: a\leq \tau\}$.

Consider the operator $D_H:A(n)\to W(n):={\rm Der}_F A(n)$ defined as
$$D_H(f)=\sum_{j=1}^{2m} \sigma(j)D_j(f)D_{j'}=\sum_{i=1}^m\left[D_i(f)D_{i+m}-D_{i+m}(f)
D_i\right],$$
where, as usual, $D_i:=\frac{\partial}{\partial x_i}\in W(n)$.

We denote with $K'(n)$ the graded Lie algebra over $F$ whose underlying $F$-vector space is 
$A(n)=F[x_1,\cdots,x_n]/(x_1^p,\cdots, x_n^p)$, 
endowed with the grading defined
by ${\rm deg}(x^a)=|a|+a_n-2$ and with the Lie bracket defined by 
\begin{equation*}\label{K-brac}
[x^a,x^b]=D_H(x^a)(x^b)+\left[a_n {\rm deg}(x^b)-b_n {\rm deg}(x^a)\right]
x^{a+b-\epsilon_n}.
\end{equation*}

\begin{defi}
The Contact algebra is the derived subalgebra of $K'(n)$:
$$K(n):=K'(n)^{(1)}=[K'(n),K'(n)].$$
\end{defi}
We need the following characterization of $K(n)$
(see \cite[Chap. 4, Theo. 5.5]{FS}).
\begin{pro}\label{K-dim}
Denote with $K'(n)_{<\tau}$ the sub-vector space of $K'(n)$ generated over $F$
by the monomials $x^a$ such that $a<\tau$. Then
 $$K(n)=
\begin{sis}
&K'(n)& \quad \text{ if } p\not | \: (m+2),\\
&K'(n)_{<\tau}& \quad \text{ if } p\: | \:(m+2).\\
\end{sis}$$
\end{pro}

We can describe explicitly the low degree terms of $K(n)$ together with their adjoint
action. The negative graded pieces of $K(n)$ are
$K(n)_{-2}=\langle 1 \rangle_F$ whose adjoint action is like the action of $2D_n$ on $A(n)$ and
$K(n)_{-1}=\oplus_{i=1}^{2m}\langle x_i \rangle_F$ where the adjoint action of $x_i$ 
is like $\sigma(i)D_{i'}+x_i D_n $. The piece $K(n)_0$ of degree $0$ is
generated by the central element $x_n$ whose adjoint action is given by 
$[x_n,x^a]={\rm deg}(x^a)x^a$
and by $x_ix_j$ (with $1\leq i,j \leq 2m$) whose adjoint action is like $\sigma(i)x_jD_{i'}+
\sigma(j)x_iD_{j'}$. 
Hence $K(n)_0\cong \mathfrak{sp}(2m, F)\oplus  \langle x_n \rangle_F$.

The algebra $K(n)$ admits a root space decomposition with respect to a canonical 
Cartan subalgebra.

\begin{pro}\label{Cartan-K}
\begin{itemize}
\item[(a)] $T_K:=\oplus_{i=1}^m \langle x_ix_{i'}\rangle_F\oplus 
\langle x_n\rangle_F$
is a maximal torus of $K(n)$ (called the canonical maximal torus).
\item[(b)] The centralizer of $T_K$ inside $K(n)$ is the subalgebra
$C_K=\{x^a \: | \: a_i=a_{i'} \: \: \text{ and } \: \: {\rm deg}(x^a)\equiv 0 \mod p\}$,
which is hence a Cartan subalgebra  (called the canonical Cartan subalgebra).
The dimension of $C_K$ is $p^m$ if $p\not | (m+2)$ and $p^m-1$ otherwise.
\item[(c)] Let $\Phi_K:={\rm Hom}_{\mathbb{F}_p}(\oplus_{i=1}^n \langle  x_ix_{i'}
\rangle_{\mathbb{F}_p}\oplus \langle x_n \rangle_{\mathbb{F}_p}, \mathbb{F}_p)$, 
where $\mathbb{F}_p$ is the prime
field of $F$. We have a Cartan decomposition $K(n)=C_K\oplus_{\phi\in \Phi_K-\0}
K(n)_{\phi}$, where $K(n)_{\phi}=\{x^a\: | \: a_{i+m}-a_i\equiv \phi(x_ix_{i'})
\: \: \forall \: i=1,\cdots, m \: \text{ and } \: {\rm deg}(x^a)\equiv\phi(x_n)\}$. 
The dimension of every $K(n)_{\phi}$, with $\phi\in \Phi_K-\0$, is $p^m$.
\end{itemize}
\end{pro}
\begin{proof}
See \cite[Chap. 4, Theo. 5.6 and 5.7]{FS}.
\end{proof}

\subsection{Proof of the Main Theorem \ref{K-finaltheorem}}

In this section, assuming the results of the next section, 
we give a proof of the Main Theorem \ref{K-finaltheorem}.

\begin{proof}[Proof of the Main Theorem \ref{K-finaltheorem}]

It is easy to see that the cochains appearing in Theorem \ref{K-finaltheorem}
are cocycles and that they are independent in $H^2(K(n), K(n))$.
Therefore, we are left with showing that $\dim_F H^2(K(n),K(n))=n$.
We divide the proof in three steps.

\underline{STEP I}: It is enough to show that $\dim_F H^2(K(n), K'(n))=n$ since
 there is an inclusion
$$H^2(K(n),K(n))\hookrightarrow H^2(K(n),K'(n)).$$
Indeed, if $p$ does not divide $m+2$ then $K'(n)=K(n)$ and we get the equality. 
Otherwise there is an exact sequence of $K(n)$-modules
\begin{equation}\label{K-redu-K'}
0\to K(n)\to K'(n)\to \langle x^{\tau} \rangle_F \to 0
\end{equation}
where $\langle x^{\tau} \rangle_F\cong F$ is the trivial $K(n)$-module. 
We get the desired inclusion since $H^1(K(n),F)=0$, which follows
from the fact that  $[K(n),K(n)]=K(n)$.

\underline{STEP II}: We have that 
$$H^2(K(n),K'(n))=H^2(K(n)_{\geq 0},F_{\lambda-\sigma}),$$
where $F_{\lambda-\sigma}$ is the one-dimensional representation of $K(n)_{\geq 0}$ on which
$x_n$ acts as $-2$ and all the others elements act trivially.

This follows from the general results of \cite{FS91}.
Indeed, it is easily seen that $K'(n)$ is the restricted $K(n)$-module induced
from the restricted $K(n)_{\geq 0}$-submodule $\langle x^{\tau}\rangle_F\subset K'(n)$.
In the notation of \cite{FS91}, 
the $K(n)_{\geq 0}$-module $\langle x^{\tau}\rangle_F$ is isomorphic to $F_{\lambda}$, 
where $F_{\lambda}$ is the one dimensional $K(n)_{\geq 0}$-module corresponding 
to the Lie algebra homomorphism $\lambda:K(n)_{\geq 0}\to F$ 
whose only non-zero value is $\lambda(x_n)=-2m-4\equiv \deg x^{\tau}=(2m+2)(p-1)-2\mod p$.

Moreover, consider the Lie algebra homomorphism $\sigma: K(n)_{\geq 0}\to F$ given by
$\sigma(x):={\rm tr}({\rm ad}_{K(n)/K(n)_{\geq 0}} x)$ for $x\in K(n)_{\geq 0}$
(see \cite[Pag. 155]{FS91}). It is easily seen that the only non-zero value of 
$\sigma$ is given by $\sigma(x_n)=-2m-2$. 

Therefore we have that $\gamma-\sigma: K(n)_{\geq 0}\to F$ 
is the Lie algebra homomorpshim sending $x_n$ to 
$-2$ and vanishing on the other elements. Moreover, using Lemma \ref{K-commutators}, 
it is straithforward to check that, in the notation of \cite{FS91}, we have  
the equality
$$\{[x, y]-(\lambda-\sigma)(x)y+(\lambda-\sigma)(y)x\: | \: x, y \in K(n)_{\geq 0}\}:=
(K(n)_{\geq 0})^{(1)}_{\lambda-\sigma}=I:={\rm ker}(\lambda-\sigma).$$
We conclude using \cite[Thm. 3.6(1)]{FS91}.

\underline{STEP III}: In section 2.3, we prove that  
$$\dim_F H^2(K(n)_{\geq 0},F_{\lambda-\sigma})=
\dim_F H^2(K(n)_{\geq 1},F_{\lambda-\sigma})^{K(n)_0}=n.$$

\end{proof}

\subsection{Computation of $H^2(K(n)_{\geq 0}, F_{\lambda-\sigma})$}

This section is devoted to complete the third step of the proof of the Main Theorem
as outlined in section 2.2, that is the computation of $H^2(K(n)_{\geq 0}, F_{\lambda-\sigma})$.
This is done in the Propositions \ref{K-H^1} and \ref{K-H^2} below.

Recall that $F_{\lambda-\sigma}$ is the one-dimensional representation of 
$K(n)_{\geq 0}$ on which $x_n$ acts as $-2$ and all the others elements act trivially.
This action becomes homogenous with respect to the weight decomposition of $K(n)_{\geq 0}$ 
if we give the weight $-2 \epsilon_n$ to the generator of $F_{\lambda-\sigma}$.
As remarked in \cite[Sec. 2.1]{Viv}, in this situation we have that 
$$H^2(K(n)_{\geq 0}, F_{\lambda-\sigma})=H^2(K(n)_{\geq 0}, F_{\lambda-\sigma})_{\0},$$
where the subscript $\0$ means that we consider only homogeneuos cochains with respect 
to the natural action of the maximal torus $T_K$ (see Prop. \ref{Cartan-K}).

\begin{pro}\label{K-H^1}
We have that 
$$ H^2(K(n)_{\geq 0},F_{\lambda-\sigma})=
H^2(K(n)_{\geq 1},F_{\lambda-\sigma})^{K(n)_0}.$$
\end{pro}
\begin{proof}
Consider the Hochschild-Serre spectral sequence (see \cite{HS}) associated to the ideal 
$K(n)_{\geq 1} \lhd K(n)_{\geq 0}$:
$$E_2^{r,s}=H^r(K(n)_0,H^s(K(n)_{\geq 1},F_{\lambda-\sigma}))
\Rightarrow H^{r+s}(K(n)_{\geq 0},F_{\lambda-\sigma}).$$
We are going to prove that the first two lines of the above 
spectral sequence vanish, which clearly imply the Proposition.

The first line $E_2^{*,0}=H^*(K(n)_0,F_{\lambda-\sigma})$ vanish for homogeneity reasons. 
Indeed, the weight of $F_{\lambda-\sigma}$ is $-2\epsilon_n\neq 0$, while the weights occurring on 
$K(n)_0$ are $\{\pm \epsilon_i \pm
\epsilon_j, \: 1\leq i,j \leq 2m\}$ and hence the weights that occur on
$K(n)_0^{\otimes k}$ cannot contain $\epsilon_n$
with a non-trivial coefficient.

On the other hand, since $F_{\lambda-\sigma}$ is a trivial $K(n)_{\geq 1}$-module and 
$[K(n)_{\geq 1},K(n)_{\geq 1}]$ $=K(n)_{\geq 2}$ by the Lemma \ref{K-commutators} below,
we have that
$$H^1(K(n)_{\geq 1},F_{\lambda-\sigma})=C^1(K(n)_1,F_{\lambda-\sigma}).$$ 
From this equality, we deduce that the second line $E_2^{*,1}=H^*(K(n)_0,H^1(W(n)_{\geq 1},$
$ F_{\lambda-\sigma}))$
vanish again for homogeneity reasons. Indeed the $n$-component of the weights appearing in 
$H^1(K(n)_{\geq 1},F_{\lambda-\sigma})=C^1(K(n)_1,F_{\lambda-\sigma})$ is $-3 \epsilon_n\neq 0$ 
(because $p\geq 5$) while the weights appearing in $K(n)_0^{\otimes k}$ have trivial $n$-component.
\end{proof}

\begin{lem}\label{K-commutators}
Let $d$ be an integer greater or equal to $-2$. Then
$$[K(n)_1,K(n)_d]=K(n)_{d+1}.$$
\end{lem}
\begin{proof}
The inclusion $[K(n)_1,K(n)_d]\subset K(n)_{d+1}$ is clear. In order to prove the other 
inclusion, we consider an element $x^a\in K(n)_{d+1}$ and we have to show that 
it belongs to the commutators $[K(n)_1,K(n)_d]$.

The elements of $K(n)_1$  are of the form $x_ix_jx_k$ or $x_ix_n$ 
(for some $1\leq i,j,k\leq 2m$). The first ones act, 
via adjoint action, as $D_H(x_ix_jx_k)-x_ix_jx_kD_n$ while the latter ones act as 
$\sigma(i)x_nD_{i' }+x_i{\rm deg}-x_ix_nD_n$.
Consider the decomposition $K'(n)=\oplus_{k=0}^{p-1}A(2m)x_n^k$. The proof is by induction
on the coefficient $a_n$, which in what follows it is called the $x_n$-degree of $x^a$.

First of all consider the case of $x_n$-degree equal to $0$,
that is the case $x^a\in A(2m)$.
If there exists an index $1\leq i\leq 2m$ such that $a_i\geq 2$ and 
$a_{i'}<p-1$, then we conclude by mean of the following formula
$$[x_i^3,x^{a-2\epsilon_i+\epsilon_{i'}}]=3\sigma(i)(a_{i'}+1)x^a.$$
Therefore it remains to consider the elements $x^a$ for which 
$a_i=a_{i'}=p-1$ or $0\leq a_i,a_{i'}\leq 1$ for every $1\leq i\leq 2m$.
If there exists a couple $(a_i,a_{i'})=(1,1)$, we are done by the formula
$$[x_i^2x_{i'},x^{a-\epsilon_i}]=\sigma(i)(2a_{i'}-a_i+1)x^a=
2 \sigma(i)x^a.$$
If there exists a couple $(a_i,a_{i'})=(1,0)$, then there are two possibilities:
either $x^a=x_i$ or there exists an index $j\neq i, i'$ such that 
$a_j\geq 1$. In the first case we use $[x_ix_n,1]=-2x_i$ while in the second we
conclude by mean of the following formula 
$$[x_i^2x_j,x^{a-\epsilon_i+\epsilon_{i'}-\epsilon_j}]=
2\sigma(i)x^a+\sigma(j)a_{j'}x^{a+\epsilon_i+\epsilon_{i'}-\epsilon_j
-\epsilon_{j'}},$$
together with the fact that the second element in the right hand side belongs to
$[K(n)_1,K(n)_d]$ by what proved above. 
Hence we are left with considering the elements $x^a$ for which every couple of 
conjugated coefficients $(a_i,a_{i'})$ is equal to $(0,0)$ or $(p-1,p-1)$. If there  
are two indexes $1\leq i\neq j\leq m$
such that $(a_i,a_{i'})=(0,0)$ and $(a_j,a_{j'})=(p-1,p-1)$ we use the formula
$$[x_j^2x_i,x^{a+\epsilon_{i'}-2\epsilon_j}]=2(p-1)x^{a+\epsilon_i
+\epsilon_{i'}-\epsilon_j-\epsilon_{j'}}+x^a,$$
together with the fact that the first term on the right hand side belongs to $[K(n)_1,$ $K(n)_d]$
by what proved above. Since the case $x^a=1$ is excluded by the hypothesis $d+1\geq -1$, 
it remains
to consider the element $x^a=x^{\sigma}$ for which we can take an appropriate linear combination 
of the two equations (with $k=0$):
\begin{equation}\label{for1-sigma}
[x_i^3,x^{\sigma-3\epsilon_i}x_n^{k+1}]=-3\sigma(i)x^{\sigma-\epsilon_i-\epsilon_{i'}}x_n^{k+1}
-(k+1)x^{\sigma}x_n^k, 
\end{equation}
\begin{equation}\label{for2-sigma}
[x_i^2x_{i'},x^{\sigma-2\epsilon_i-\epsilon_{i'}}x_n^{k+1}]=-\sigma(i)x^{\sigma-\epsilon_i-
\epsilon_{i'}}x_n^{k+1}-(k+1) x^{\sigma}x_n^k.
\end{equation}
For the inductive step, suppose that $a_n=k\geq 1$ and that we have already proved the 
desired inclusion for the elements of $x_n$-degree less than or equal to $k-1$.
If there exists
an index $i$ such that $a_i<p-1$, then the formula
$$[x_{i'}x_n,x^{a+\epsilon_{i}-\epsilon_n}]=\sigma(i')(a_i+1)x^a+(d-a_n+1)x^{a+\epsilon_i
+\epsilon_{i'}-\epsilon_n},$$
together with the induction hypothesis, gives the conclusion. Otherwise our element is 
equal to $x^{\sigma}x_n^k$. If $k<p-1$, then one concludes by taking an appropriate 
linear combination of the above formulas (\ref{for1-sigma}) and (\ref{for2-sigma}).
Finally for the element $x^{\sigma}x_n^{p-1}=x^{\tau}$ (which can occur only if 
$p\not | m+2$), the conclusion follows from the formula (for an arbitrary chosen
$1\leq i\leq 2m$)
$$[x_ix_n,x^{\sigma-\epsilon_i}x_n^{p-1}]=-2(m+2)x^{\tau}.$$ 
\end{proof}

\noindent For the remaining part of this subsection, we identify the $K(n)_0$-module 
$F_{\lambda-\sigma}$ with 
the $K(n)_0$-module $F\cong\langle 1\rangle_F=K(n)_{-2}$. On both these modules, 
$K(n)_{\geq 1}$ acts trivially.

\begin{pro}\label{K-H^2}
We have that
$$H^2(K(n)_{\geq 1},F)^{K(n)_0}=\bigoplus_{i=1}^{2m} 
\langle \overline {{\rm Sq}(x_i)} \rangle_F \oplus \langle \overline{\Sq(1)} \rangle_F, $$ 
where $\overline{{\rm Sq}(x_i)}$ is the projection of ${\rm Sq}(x_i)$ onto 
$\langle 1\rangle_F\cong F $ (analogously for $\overline{{\rm Sq}(1)}$).
\end{pro}
\begin{proof} 
It is easy to check that the above cocycles are independent modulo coboundaries, 
so we have to prove that they generate the whole cohomology group.

The strategy of the proof is exactly the same as that of proposition 
\cite[Prop. $3.10$]{Viv}, that is to compute, 
step by step as $d$ increases, the truncated invariant cohomology groups
$$H^2\left(\frac{K(n)_{\geq 1}}{K(n)_{\geq d+1}},F\right)^{K(n)_0}.$$
Observe that if $d$ is big enough (at least $2(m+1)(p-1)-1)$ then $K(n)_{\geq d+1}=0$ and
hence we get the cohomology we are interested in. On the other hand, by homogeneity, we 
get that   
$$H^2\left(\frac{K(n)_{\geq 1}}{K(n)_{\geq 2}},F\right)^{K(n)_0}=
C^2\left(K(n)_1,F\right)^{K(n)_0}=0. $$
By taking the Hochschild-Serre spectral sequence (see \cite{HS}) associated to the ideal 
$K(n)_d=\frac{K(n)_{\geq d}}{K(n)_{\geq d+1}}\lhd \frac{K(n)_{\geq 1}}{K(n)_{\geq d+1}}$
(for $d\geq 2$):
\begin{equation}\label{K-HS-final}
E_2^{r,s}=H^r\left(\frac{K(n)_{\geq 1}}{K(n)_{\geq d}},H^s(K(n)_d,F)\right)
\Rightarrow H^{r+s}\left(\frac{K(n)_{\geq 1}}{K(n)_{\geq d+1}},F\right),
\end{equation}
we get the same diagram as in \cite[Prop. 3.10]{Viv}
 (the vanishing of $E_2^{0,2}$ and the injectivity of the map $\alpha$ are proved in exactly 
the same way).

By taking the cohomology with respect to $K(n)_0$ and using the Lemmas \ref{K-inv2-(1,1)},
\ref{K-inv-(0,1)}, \ref{K-H^1bis-(0,1)} below, we see that the only cocycles that contribute
to the required cohomology group are 
$\{\overline {{\rm Sq}(x_1)},\cdots, \overline{{\rm Sq}(x_{2m})},
\overline{{\rm Sq}(1)}\}$ since the cocycle ${\rm inv}\circ [-,-]\in 
(E_{\infty}^{1,1})^{K(n)_0}$ in degree $2m(p-1)+2\mu-3$ is annihilated by ${\rm inv}\in 
C^1(K(n)_{2m(p-1)+2\mu-2}, F)^{K(n)_0}$. 
\end{proof}

The remaining part of this section is devoted to the proof of the Lemmas that were used 
in the proof of Proposition \ref{K-H^2}. 

In the next two Lemmas, we are going to compute the $K(n)_0$-invariant terms 
$(E_2^{1,1})^{K(n)_0}$ and  $(E_{\infty}^{1,1})^{K(n)_0}$. Observe that,
since $K(n)_d$ is in the center of $\frac{K(n)_{\geq 1}}{K(n)_{\geq d+1}}$,
we have that $H^s(K(n)_d,F)=C^s(K(n)_d,F)$ and $\frac{K(n)_{\geq 1}}{K(n)_{\geq d+1}}$
acts trivially on it. Therefore, using Lemma \ref{K-commutators}, we deduce that 
$E_2^{1,1}=C^1(K(n)_1\times K(n)_d,F)$.

\begin{lem}\label{K-inv-(1,1)}
Let $0\leq \mu,\nu\leq p-1$ such that $\mu\equiv m \mod p$ and $\nu\equiv (m+1) \mod p$.
Then we have that
$$\bigoplus_{d\geq 2} C^1(K(n)_1\times K(n)_d,F)^{K(n)_0}=\langle \Phi_1\rangle_F \oplus 
\langle \Phi_2\rangle_F \oplus \langle \Psi_1\rangle_F \oplus \langle \Psi_2\rangle_F,$$
where the above cochains are defined by 
$$\begin{sis}
&\Phi_1(x_ix_n, x_{i'}x_n^{p-1})=\sigma(i),\\
&\Phi_2(x^{a}, x^{\widehat{a}}x_n^{p-2})=\sigma(a)a! \hspace{1cm} \text{ if } \: a\in \N^{2m}
\: \text{ and } |a|=3, \\
&\Psi_1(x_ix_n, x^{\sigma-\epsilon_i}x_n^{\mu})=1,\\
&\Psi_2(x^a,x^{\sigma-a}x_n^{\nu})=1 \hspace{1,7cm} \text{ if } \: a\in \N^{2m}
\: \text{ and } |a|=3. 
\end{sis}$$
\end{lem}
\begin{proof}
An easy verification shows that the four cochains of above are $K(n)_0$-invariants
and linearly independent. 
We will conclude by showing that the dimension over the base field $F$ of the space 
$\bigoplus_{d\geq 2} C^1\left(K(n)_1\times  K(n)_d, F\right)_{\0}^{K(n)_0}$ of all 
invariant homogeneous cochains is less than or equal to $4$.

The space $K(n)_1$ 
admits the decomposition $K(n)_1=A(2m)_{-1}\cdot x_n\oplus A(2m)_1$ which is invariant
under the adjoint action of $K(n)_0=A(2m)_0\oplus \langle x_n \rangle_F$. Moreover 
the action of $K(n)_0$ is transitive in both the summands 
$A(2m)_{-1}\cdot x_n$ and $A(2m)_1$. Therefore a $K(n)_0$-invariant homogeneous cochain 
$g\in  C^1\left(K(n)_1\times  K(n)_{\geq 2},F \right)_{\0}^{K(n)_0}$ 
is determined by the values on any two elements of $A(2m)_{-1}\cdot x_n$ and $A(2m)_1$, 
let's say $x_1x_n$ and $x_1^3$.   

Consider an element $x^a\in K(n)_{\geq 2}$ such that $g(x_1x_n,x^a)\neq 0$.
By homogeneity the element $x^a$ must satisfy $a_{1'}\equiv a_1+1 \mod p$,
$a_{j'}=a_j$ for every $j\not\in \{1,1'\}$ and ${\rm deg}(x^a)\equiv -3 \mod p$. 
If the couple $(a_1,a_{1'})$ would be different 
from $(0,1)$ or $(p-2,p-1)$ then the following invariance condition
$$0=(x_1^2\circ g)(x_1x_n, x^{a-\epsilon_1+\epsilon_{1'}})=-2(a_{1'}+1)g(x_1x_n,x^a)$$
would contradict the hypothesis of non-vanishing. Therefore we can assume that 
$(a_1,a_{1'})=(0,1)$ or $(p-2,p-1)$. If the first case holds, then necessarily 
$(a_j,a_{j'})=(0,0)$ for every $j\not\in\{1,1'\}$. Indeed if this is not the case, then
we get a contradiction with the non-vanishing hypothesis by means of the following 
invariance condition
$$0=(x_1x_j\circ g)(x_1x_n,x^{a+\epsilon_{1'}-\epsilon_j})=-2g(x_1x_n, x^a).$$ 
Analogously, if $(a_1,a_{1'})=(p-2,p-1)$ then $(a_j,a_{j'})$ for every $j\not\in\{1,1'\}$
because  of the following invariance condition
$$0=(x_1x_j\circ g)(x_1x_n,x^{a-\epsilon_1+\epsilon_{j'}})=-\sigma(j)(a_ {j'}+1)
g(x_1x_n,x^a).$$ 
Taking into account the homogeneity condition ${\rm deg}(x^a)\equiv -3 \mod p$, we get 
that the only non-zero values of $g(x_1x_n,-)$ can be $g(x_1x_n,x_{1'}x_n^{p-1})$ and
$g(x_1x_n,x^{\sigma-\epsilon_1}x_n^{\mu})$.

In exactly the same way, one proves that the only non-zero values of $g(x_1^3,-)$ 
can be $g(x_1^3, x_{1'}^3x_n^{p-2})$ and $g(x_1^3, x^{\sigma-3\epsilon_1}x_n^{\nu})$ and 
therefore we get 
$${\rm dim}_F\,  C^1\left(K(n)_1\times  K(n)_{\geq 2},F\right)_{\0}^{K(n)_0}\leq 4.
$$
\end{proof}

\begin{lem}\label{K-inv2-(1,1)}
In the above spectral sequence (\ref{K-HS-final}), we have that
$$(E_{\infty}^{1,1})^{K(n)_0}=
\begin{cases}
\langle\overline{{\rm Sq}(1)}\rangle_F & \text{ if } d=2p-3, \\
\langle {\rm inv}\circ [-,-]\rangle_F &\text{ if } d=2m(p-1)+2\mu-3,\\
0 & \text{ otherwise, }
\end{cases}$$
where ${\rm inv}\in C^1(K(n)_{2m(p-1)+2\mu-2},F)^{K(n)_0}$ is defined in the
Lemma \ref{K-inv-(0,1)} below and $\overline{{\rm Sq}(1)}$ is the restriction of 
${\rm Sq}(1)$ to $K(n)_1\times K(n)_{2p-3}$.
\end{lem}
\begin{proof}
$(E_{\infty}^{1,1})^{K(n)_0}$ is the subspace of $(E_2^{1,1})^{K(n)_0}=
C^1(K(n)_1\times K(n)_d,F)^{K(n)_0}$
consisting of cocycles that can be lifted to $Z^2\left(\frac{K(n)_{\geq 1}}{K(n)_{\geq d+1}},
F\right)$.
A direct computation shows that, with the notations of Lemma \ref{K-inv-(1,1)}, 
$\overline{{\rm Sq}(1)}=2\overline{{\rm Sq}(D_n)}=2 \Phi_1$ and ${\rm inv}\circ [-,-]=
-(\nu+2)\Psi_1-\nu\Psi_2$ and clearly these
two cocycles can be lifted. We want two show that any other liftable cocycle is a linear 
combination of them.

First of all we show that the cocycle $\Phi_2$ is not liftable. By absurd, 
suppose that a lifting exists and call it again $\Phi_2$. We get a contradiction by 
mean of the following cocycle conditions
$$\begin{sis}
& 0=d\Phi_2(x_i^3,x_{i'}^3,x_n^{p-1})=-9\sigma(i)\Phi_2(x_i^2x_{i'}^2,x_n^{p-1})+
\Phi_2(x_i^3x_n^{p-2},x_{i'}^3)+\\
& -\Phi_2(x_{i'}^3x_n^{p-2},x_i^3)=-9\sigma(i)\Phi_2(x_i^2x_{i'}^2,x_n^{p-1})+12\sigma(i),\\
& 0= d\Phi_2(x_i^2x_{i'},x_{i'}^2x_i,x_n^{p-1})=-3\sigma(i)\Phi_2(x_i^2x_{i'}^2,x_n^{p-1})+
\Phi_2(x_i^2x_{i'}x_n^{p-2},x_{i'}^2x_i)+\\
&-\Phi_2(x_{i'}^2x_ix_n^{p-2},x_i^2x_{i'})= -3\sigma(i)\Phi_2(x_i^2x_{i'}^2,x_n^{p-1})
-4\sigma(i).
\end{sis}$$
Finally, consider the cocycle $\Psi\in\langle \Psi_1, \Psi_2\rangle_F$ defined in degree $2m(p-1)+
2\nu-5$ by 
$$\Psi=
\begin{cases}
a\Psi_1+b\Psi_2 & \text{ if } \nu\neq 0,\\
b \Psi_2 & \text{ if } \nu=0,
\end{cases}$$ 
for certain $a,b\in F$. 
We will show that $\Psi$ can be lifted to $Z^2\left(\frac{K(n)_{\geq 1}}
{K(n)_{\geq 2m(p-1)+2\nu-4}},F
\right)$ if and only if $b(\nu+2)\equiv a\nu \mod p$ and this will conclude our proof. 
Indeed this 
imply that if $\nu\neq 0$ then $\Psi$ is liftable if and only if it is a multiple of ${\rm inv}
\circ[-,-]$, while if $\nu=0$ it implies that $\Psi_2$ is not liftable and hence again that 
${\rm inv}\circ[-,-]=-2 \Psi_1$ is the only liftable cocycle in the span of $\Psi_1$ and 
$\Psi_2$.

So suppose that a lift exists and call it again $\Psi$. From the following cocycle condition
$$\begin{aligned}
& 0=d\Psi(x_i^3,x_{i'}^3,x^{\sigma-2\epsilon_{i'}-2\epsilon_i}x_n^{\nu})=
-9\sigma(i)\Psi(x_i^2x_{i'}^2,x^{\sigma-2\epsilon_{i'}-2\epsilon_i}x_n^{\nu})+\\
&-9\sigma(i)\Psi(x^{\sigma-3\epsilon_{i'}}x_n^{\nu},x_{i'}^3)+
9\sigma(i')\Psi(x^{\sigma-3\epsilon_i}x_n^{\nu},x_i^3),
\end{aligned}$$
we deduce that $\Psi(x_i^2x_{i'}^2,x^{\sigma-2\epsilon_{i'}-2\epsilon_i}x_n^{\nu})=2b$. 
Using this, we get the following 
$$\begin{aligned}
&0=d\Psi(x_ix_n,x_ix_{i'}^2, x^{\sigma-2\epsilon_i-2\epsilon_{i'}}x_n^{\nu})=
-2\sigma(i)\Psi(x_ix_{i'}x_n+x_i^2x_{i'}^2,x^{\sigma-2\epsilon_i-2\epsilon_{i'}}x_n^{\nu})+\\
&+(-4-\nu)\Psi(x^{\sigma-\epsilon_i-2\epsilon_{i'}}x_n^{\nu},x_ix_{i'}^2)+
\nu\Psi(x^{\sigma-\epsilon_i}x_n^{\mu},x_ix_n)=\\
&=-2\sigma(i)\Psi(x_ix_{i'}x_n,x^{\sigma-2\epsilon_i-2\epsilon_{i'}}x_n^{\nu})+(\nu+2)b-\nu a.
\end{aligned}$$
Exchanging $i$ with $i'$ and summing the two expressions, we obtain the required 
congruence $(\nu+2)b\equiv \nu a\mod p$. 
\end{proof}

In the next Lemma, we compute the $K(n)_0$-invariants of the term $E_2^{0,1}=
C^1(K(n)_d,F)$ of the above spectral sequence (\ref{K-HS-final}).

\begin{lem}\label{K-inv-(0,1)}
Let $\mu$ the integer defined in Lemma \ref{K-inv-(1,1)}.
We have that 
$$C^1\left(K(n)_d, F\right)^{K(n)_0}=
\begin{cases}
\langle {\rm inv} \rangle_F & \text{ if } d=2m(p-1)+2\mu-2,\\
0& \text{ otherwise, }
\end{cases}$$
where ${\rm inv}\in C^1(K(n)_{2m(p-1)+2\mu-2},F)$ sends $x^{\sigma}x_n^{\mu}$
into $1$ and vanish on the other elements.
\end{lem}
\begin{proof}
First of all observe that if $p$ divide $(m+2)$, then $\mu=p-2$ and hence 
$x^{\sigma}x_n^{\mu}\in K(n)$. The (well-defined) cochain
${\rm inv}$ is $K(n)_0$-invariant.
Indeed it is homogeneous and the invariance with respect to an element $x_ix_j\in K(n)_0
\setminus T_K$
(hence with $j\neq i'$) follows from the fact that $[x_ix_j,x^ {\sigma}x_n^{\mu}]=[x_ix_j,1]=0$ 
together with the fact that $x^{\sigma}x_n^{\mu}\not\in [x_ix_j,K(n)]$.  

Consider next an invariant cochain $f\in C^1\left(K(n)_d, F\right)^{K(n)_0}$ and let $x^a\in
K(n)_d$ be an element such that $f(x^a)\neq 0$. Then by homogeneity it must hold
that $a_i=a_{i'}$ for every $1\leq i\leq m$ and ${\rm deg}(x^a)\equiv -2 \mod p$.
Using the invariance with respect to $x_i^2$ or $x_{i'}^2$, we obtain that $a_i=a_{i'}=0$ 
or $p-1$. Otherwise, assuming, up to interchanging $i$ with $i'$,
that $a_i>0$ and $a_{i'}<p-1$, one gets the vanishing as follows 
$$0= (x_i^2\circ f)(x^{a-\epsilon_i+\epsilon_{i'}})=-2\sigma(i)(a_{i'}+1)f(x^a).$$
Moreover, if there are two pairs verifying $(a_i,a_{i'})=(0,0)$ and $(a_j,a_{j'})=
(p-1,p-1)$ (for $j\neq i,i'$), then we obtain the vanishing 
by means of the following 
$$0=(x_ix_j\circ f)(x^{a+\epsilon_{i'}-\epsilon_j})=-\sigma(i)f(x^a).$$
Finally, by imposing ${\rm deg}(x^a)\equiv -2\mod p$, we deduce that $x^a=1$ (which we can 
exclude since ${\rm deg}(x^a)=d\geq 2$) or $x^a=x^{\sigma}x_n^{\mu}$.
\end{proof}

In the next Lemma, we compute the first cohomology group with respect to $K(n)_0$
of the term $E_2^{0,1}=C^1(K(n)_d,F)$ of the above spectral sequence (\ref{K-HS-final}).

\begin{lem}\label{K-H^1-(0,1)}
Let $\mu$ be the integer defined in Lemma \ref{K-inv-(1,1)}. We have that
$$H^1\left(K(n)_0,C^1\left(K(n)_d,F\right)\right)=
\begin{cases}
\bigoplus_{i=1}^{2m} \langle \overline{{\rm Sq}(x_i)} \rangle_F  \hspace{2,5cm}
\text{ if } d=p-2,\\
\bigoplus_{j=1}^{2m} \langle \omega_j \rangle_F  \hspace{0,4cm}
\text{ if } d=2m(p-1)+2\mu-2-p,\\
\langle x_n \mapsto {\rm inv} \rangle_F \hspace{0,8cm}\text{ if } d=2m(p-1)+2\mu-2,\\
0  \hspace{5,1cm} \text{ otherwise, }
\end{cases}$$
where $\overline{{\rm Sq}(x_i)}$ denotes the restriction of ${\rm Sq}(x_i)$ to 
$K(n)_0\times K(n)_{p-2}$ and the cocycle $\omega_i$ is defined by 
(with $j\neq i, i'$)
$$\begin{sis}
&\omega_i(x_i^2, x^{\sigma-\epsilon_i-(p-1)\epsilon_{i'}}x_n^{\mu})= 2,\\
&\omega_i(x_ix_j, x^{\sigma-(p-1)\epsilon_{i'}-\epsilon_j} x_n^{\mu})=1.
\end{sis}$$
\end{lem}
\begin{proof}
By homogeneity we can restrict to the case $d\equiv -2 \mod p$.
First of all I claim that $f_{x_ix_{i'}}=0$ for every $1\leq i\leq n$ and $f_{x_n}$ takes 
a non zero-value only on the element $x^{\sigma}x_n^{\mu}$. Indeed, by 
the homogeneity assumption, we get for an element $\gamma\in T_K$ that
$0=df_{|\gamma}=\gamma\circ f- d(f_{|\gamma})=-d(f_{|\gamma})$, that is $f_{\gamma}\in 
C^1(K(n)_d, F)^{K(n)_0}=\langle {\rm inv} \rangle_F$ (see Lemma 
\ref{K-inv-(0,1)}). 
Moreover the cocycle condition 
$$0=df_{(x_i^2, x_{i'}^2)}(x^{\sigma}x_n^{\mu})=-4\sigma(i)f_{x_ix_{i'}}(x^{\sigma}x_n^{\mu})$$
gives $f_{x_ix_{i'}}=0$ for every $1\leq i\leq n$, while the fact that $x_n\not\in 
[K(n)_0,K(n)_0]$ implies that $f_{x_n}(x^{\sigma}x_n^{\mu})$ can be different from zero.

Now we split the proof into two parts according to the cases $d=(2r+1)p-2$ or $d=2rp-2$ for some 
integer $r$.

$\underline{\text{ I CASE}: d=(2r+1)p-2.}$

Note that in this case, there are not coboundary elements since, by reasons of parity,
$C^1(K(n)_d,F)_{\0}=0$. Moreover, for a homogenoeus coycle 
$f\in C^1(K(n)_0, $ $C^1(K(n)_d,F))_{\0}$, the value $f_{x_i^2}(x^a)$ can be different 
from $0$ only if one of the following possibilities occur
\begin{align*}
& (a_i,a_{i'})=(p-1,1) \:\text{ and }\: a_j=a_{j'} \:\text{ for every } 
j\neq i,i', \tag*{(A)}\\
& (a_i,a_{i'})=(p-2,0) \:\text{ and }\: a_j=a_{j'} \:\text{ for every } j\neq i,i'.
 \tag*{(B)}
\end{align*}
Analogously, if $j\neq i, i'$, then $f_{x_ix_j}(x^a)$ can be different from $0$ only if
(up to interchanging $i$ and $j$)
\begin{equation*}
(a_i,a_{i'})=(p-1,0), \: 1\leq a_{j'}=a_j+1 \leq p-1 \: \text{ and } 
\:a_k=a_{k'} \:\text{ for }k\neq i,i'\,j,j'. \tag*{(C)}
\end{equation*}  
The values of types $(C)$ are determined by the values of types $(A)$ and $(B)$ 
by mean of the following 
cochain condition (where $a$ is a multindex as in $(C)$)
$$0=df_{(x_ix_j,x_i^2)}(x^{a-\epsilon_i+\epsilon_{i'}})=$$
$$=-\sigma(j)a_{j'}
f_{x_i^2}(x^{a+\epsilon_{i'}-\epsilon_{j'}})-\sigma(i)f_{x_i^2}(x^{a-\epsilon_i+
\epsilon_j})+ 2\sigma(i)f_{x_ix_j}(x^a)$$
where in the last equation the first term is of type $(A)$ (or vanish) and the second is 
of type $(B)$ (or vanish).

The values of type $(A)$ vanish if there exists an index $j\neq i,i'$ 
such that $a_j=a_{j'}\neq 0$, because of the following condition 
(where $a$ satisfies the conditions in $(A)$)
$$0=df_{(x_ix_j,x_i^2)}(x^{a+\epsilon_{i'}-\epsilon_j})=-2 \sigma(i) f_{x_i^2}(x^a).$$
On the other hand, the values of type $(B)$ vanish if there exists a $j\neq i,i'$ 
such that $a_j=a_{j'}\neq p-1$ because
the following cocycle condition (where $a$ satisfies the conditions of $(B)$)
$$0=df_{(x_ix_j,x_i^2)}(x^{a-\epsilon_i+\epsilon_{j'}})=-\sigma(j)(a_{j'}+1)f_{x_i^2}(x^a).$$
Therefore a cochain $f$ is completely determined by the values $f_{x_i^2}(x_i^{p-1}x_{i'})$ 
and $f_{x_i^2}(x^{\sigma-\epsilon_i-(p-1)\epsilon_{i'}}x_n^{\mu})$, whose values determine also 
the cocycles $\overline{{\rm Sq}(x_{i'})}$ and $\omega_i$ (respectively), and hence $f$ 
is a linear combination of $\overline{{\rm Sq}(x_{i'})}$ or $\omega_i$.

$\underline{\text{ II CASE}: d=2rp-2.}$

In this case we will prove that $f$ vanish (up to adding a coboundary $dg$) except for the 
value $f_{x_n}(x^{\sigma}x_n^{\mu})$
(which can be non-zero as seen before). We have already seen that $f_{x_ix_{i'}}$ vanish for 
every $1\leq i\leq n$. 

We first prove that, by adding coboundaries, we can modify the cochain $f$ (without changing
its cohomological class) in such a way that it satisfies $f_{x_i^2}=0$  
for every $1\leq i\leq m$. The proof is by induction on $i$. So suppose that for 
a certain  $k$, we have that $f_{x_i^2}=0$ for every $i<k$.
We want to prove that, by adding coboundaries, 
we can modify $f$ in such a way that it verifies $f_{x_k^2}=0$.

First of all note that, by homogeneity and parity condition on $d$, $f_{x_k^2}(x^a)$ can be 
different from $0$ only if $2\leq a_{k'}=a_k+2\leq p-1$ and $a_{h'}=a_h$ for $h\neq k,k'$.
Moreover if there exists an index $1\leq h<k\leq m$ such that $a_h=a_{h'}\neq 0, (p-1)$, then 
$f_{x_k^2}(x^a)=0$ because of the following cocycle condition
$$0=df_{(x_h^2,x_k^2)}(x^{a-\epsilon_h+\epsilon_{h'}})=-2(a_{h'}+1)f_{x_k^2}(x^a)$$
where we used that $f_{x_h^2}=0$ by induction. Therefore we can suppose that
for $1\leq h<k\leq m$, $a_h=a_{h'}=0$ or $(p-1)$. Fix one of these elements $x^a$.
Define an element 
$g\in C^1(K(n)_d,F)_{\0}$ as follows:
$$\begin{sis}
& g(x^{a+\epsilon_k-\epsilon_{k'}})=\frac{f_{x_k^2}(x^a)}{2 a_{k'}}, & \\
& g(x^b)=0 & \text{ if } b\neq a+\epsilon_k-\epsilon_{k'}.\\
\end{sis}$$
By construction, if $1\leq h< k \leq m$ then $(x_h^2\circ g)=0$ while 
$(x_k^2\circ g)(x^a)=-f_{x_k^2}(x^a)$. Therefore the new cocycle 
$\widetilde{f}:=f+dg$ satisfies the same inductive hypothesis as before and moreover
it verifies $\widetilde{f}_{x_k^2}(x^a)=0$. Repeating these modifications for all the elements 
$x^a$ as before, eventually we obtain a new cochain homologous to the old one (which,
by an abuse of notation, we continue to call $f$) and which satisfies $f_{x_i^2}=0$ for
every $1\leq i\leq k$, as required.

Using the above conditions, we want to show that the cochain $f$ must satisfy also 
$f_{x_{i'}^2}=0$ for every 
$1\leq i\leq m$ (and hence that $f_{x_j^2}=0$ for every $j$). Indeed, as before,
we have that $f_{x_{i'}^2}(x^a)$ can be different from $0$ only if 
$2\leq a_i=a_{i'}+2\leq p-1$ and $a_j=a_{j'}$ for every $j\neq i,i'$. Hence the required 
vanishing follows from the following cocycle condition 
$$0=df_{(x_i^2,x_{i'}^2)}(x^{a-\epsilon_i+\epsilon_{i'}})=-2(a_{i'}+1)f_{x_{i'}^2}(x^a).$$ 
Finally we have to show that we can modify once more (by adding coboundaries) 
the cocycle $f$ in such a way
that the previous vanishings $f_{x_i^2}=0$ are still satisfied and moreover also $f_{x_ix_j}$ 
vanish for every $j\neq i,i'$. 

First of all, note that using cocycle 
conditions of type $0=df_{(x_h^2,x_ix_j)}$ with $h\neq i', j'$ and the fact that $f_{x_h^2}=0$, 
we obtain the vanishing of $f_{x_ix_j}(x^a)$ for all the elements 
$x^a\in [x_h^2,K(n)_d]\cap K(n)_{-\epsilon_i-\epsilon_j}$ (for $h\neq i',j'$), 
that is for all the elements of $x^a\in K(n)_d$ with the exception of the ones that verify 
$$(a_k,a_{k'})=\begin{sis}
& (0,1) \text{ or } (p-2,p-1) & \text{ if } k=i \text{ or } j,\\
& (0,0) \text{ or } (p-1,p-1) & \text{ otherwise. }
\end{sis}$$
Therefore, we can assume that our $x^a$ verifies these conditions.  For the rest of the proof, 
we introduce the following definitions. We say that a couple $(a_k,a_{k'})$ is small if it is 
equal to $(0,0)$ or $(0,1)$ or $(1,0)$ according to the conditions above, while we say 
that it is big if
it is equal to $(p-1,p-1)$ or $(p-2,p-1)$ or $(p-1,p-2)$. Moreover we say 
that $x^a$ has an ascending jump in position $k$ (with $1\leq k \leq m-1$) if $(a_k,a_{k'})$
is small and $(a_{k+1},a_{(k+1)'})$  is big, while we say that it has a descending jump 
in position $k$ if $(a_k,a_{k'})$ is big and $(a_{k+1},a_{(k+1)'})$ is small.  

We want to modify our cocycle $f$, by adding coboundaries, in such a way that 
$f_{x_ix_j}(x^a)$ vanish if $x^a$ has a jump. 

We prove this for the elements $f_{x_ix_{i+1}}$ with 
$1\leq i\leq m-1$. It is enough to prove that $f_{x_ix_j}(x^a)=0$ if there is a jump in a 
position less than or equal to $i$. Indeed if the jump on $x^a$ occurs for $h>i$, then one 
obtains the vanishing using the cocycle condition $0=df_{(x_ix_{i+1},x_hx_{h+1})}$. 
Hence, by induction on $i$, suppose that we have already proved this for the elements
$i\leq k-1$ and we want to prove it for $f_{x_kx_{k+1}}$. If there is a jump in the element 
$x^a$ occurring in a position $h<k$ then the vanishing follows from a cocycle condition of type 
$0=df_{(x_hx_{h+1},x_kx_{k+1})}$ plus the induction hypothesis. If the first jump occurring
in $x^a$ is in the $k$-th position, then we define an element $g\in C^1(K(n)_d,
F)_{\0}$ as follows:
$$\begin{sis}
&g(x^{a-\epsilon_{k'}+\epsilon_{k+1}})=f_{x_kx_{k+1}}(x^a)& \text{ if the jump is 
ascending,}\\
&g(x^{a-\epsilon_{(k+1)'}+\epsilon_k})=f_{x_kx_{k+1}}(x^a)& \text{ if the jump is 
descending,}\\
&g(x^b)=0 & \text{ otherwise.}
\end{sis}$$  
By construction (and the hypothesis on $x^a$), for every $1\leq j\leq 2m$ we have that 
$(x_j^2\circ g)=0$ and if $1\leq h<k\leq m-1$ then 
$(x_hx_{h+1}\circ g)=0$ while $(x_kx_{k+1}\circ g)(x^a)=-f_{x_kx_{k+1}}(x^a)$. Therefore
the new cocycle $\widetilde{f}=f+dg$ satisfies the same vanishing conditions of $f$ (namely 
$\widetilde{f}_{x_j^2}=0$ for every $j$ and $\widetilde{f}_{x_hx_{h+1}}=0$ for $1\leq h<k$)  plus
the new one $\widetilde{f}_{x_kx_{k+1}}(x^a)=0$. Repeating these modifications for all the 
elements $x^a$ as above, we find a new cocycle (which, by an abuse of notation, we will 
still call $f$) that satisfies $f_{x_kx_{k+1}}=0$, concluding thus the inductive step.

From the previous special cases, it follows also the vanishing of  $f_{x_ix_j}(x^a)$ 
(always under the presence of a jump) if $1\leq i,j\leq m$. Indeed, if an element $x^a$ as 
before has a jump in position $k$ then the coboundary condition
$$0=df_{(x_ix_j,x_kx_{k+1})}(x^{a+\epsilon_{k'}-\epsilon_{k+1}})=-\sigma(k)(a_{k'}+1)
f_{x_ix_j}(x^a),$$ 
in the case of an ascending jump, and 
$$0=df_{(x_ix_j,x_kx_{k+1})}(x^{a+\epsilon_{(k+1)'}-\epsilon_k})=-\sigma(k+1)(a_{(k+1)'}+1)
f_{x_ix_j}(x^a),$$
in the case of a descending jump, gives the required vanishing. 

Finally, the general case 
(in which $i$ and $j$ can vary from $1$ to $2m$) follows from cocycle conditions of type 
$0=df_{(x_ix_j,x_{i'}^2)}=-(x_{i'}^2\circ f_{x_ix_j})-2\sigma(i)f_{x_{i'}x_j}$.

So it remains to consider only the elements $x^a$ without jumps or, in other words, it remains 
to prove the vanishing of the following values of $f$: 
$f_{x_ix_j}(x_{i'}x_{j'}x_n^{p-1})=\alpha_{ij}\cdot 1$ and 
$f_{x_ix_j}(x^{\sigma-\epsilon_i-\epsilon_j}x_n^{\nu})=\beta_{ij}\cdot 1$, where 
$\nu\equiv m+1\mod p$ and $0\leq\nu\leq p-1$.  The first ones vanish because of the
following two cocycle conditions
$$\begin{sis}
&0=df_{(x_{i'}^2,x_ix_j)}(x_ix_{j'}x_n^{p-1})=-2\sigma(i')\alpha_{ij}-2\sigma(i')\alpha_{i'j},\\
&0=df_{(x_ix_j,x_{i'}x_j)}(x_{j'}^2x_n^{p-1})=-2\sigma(j)\alpha_{i'j}+2\sigma(j)\alpha_{ij}.
\end{sis}$$
The second ones vanish because of the following two cocycle conditions
 $$\begin{sis}
 &0=df_{(x_{i'}^2,x_ix_j)}(x^{\sigma-\epsilon_{i'}-\epsilon_j}x_n^{\nu})=
 -2\sigma(i')(p-1)\beta_{ij}-2\sigma(i')\beta_{i'j},\\
 &0=df_{(x_ix_j,x_{i'}x_j)}(x^{\sigma-2\epsilon_j}x_n^{\nu})=-\sigma(i)(p-1)\beta_{i'j}+
 \sigma(i')(p-1)\beta_{ij}.
 \end{sis}$$    
\end{proof}

In the next (and last) Lemma, we consider the differential map 
\begin{equation}\label{K-diff-map}
\d: E_2^{0,1}=C^1(K(n)_d,F)\to E_2^{2,0}=H^2\left( \frac{K(n)_{\geq 1}}{K(n)_{\geq d}},F
\right),
\end{equation}
induced by the above spectral sequence (\ref{K-HS-final}).
We compute the kernel of the induced map on the first cohomology group with respect to
$K(n)_0$.

\begin{lem}\label{K-H^1bis-(0,1)}
Consider the map 
$${\rm d}^{(1)}:H^1(K(n)_0,C^1(K(n)_d,F))\longrightarrow H^1\left(K(n)_0,H^2\left
(\frac{K(n)_{\geq 1}}{K(n)_{\geq d}},F\right)\right)$$
induced by the differential map (\ref{K-diff-map}). The kernel of $\d^{(1)}$ is given by
$${\rm Ker}({\rm d})=\begin{cases}
\bigoplus_{i=1}^{2m} \langle \overline{{\rm Sq}(x_i)}\rangle_F 
& \text{ if } d=p-2, \\
0 & \text{ otherwise, }
\end{cases}$$ 
where $\overline{{\rm Sq}(x_i)}$ denotes the restriction of ${\rm Sq}(x_i)$ to $K(n)_0\times 
K(n)_{p-2}$.
\end{lem}
\begin{proof}
Clearly the cocycles $\overline{{\rm Sq}(x_i)}$, being the restriction of global cocycles, 
belong to the kernel of ${\rm d}$. 
We want to show that the other generators of $H^1(K(n)_0,$ $C^1(K(n)_d,F))$ (see Lemma 
\ref{K-H^1-(0,1)}) does not belong to ${\rm Ker}({\rm d}^{(1)})$.
First of all we have that
$${\rm d}^{(1)}\langle x_n\mapsto {\rm inv}\rangle_F =\langle x_n \mapsto 
{\rm inv}\circ [-,-]\rangle,$$
and this last cocycle is not a coboundary 
since ${\rm inv}\circ [-,-]\in H^2\left(\frac{K(n)_{\geq 1}}{K(n)_{\geq d}},F\right)_{\0}$ and 
$x_n$ acts trivially on this space.

Consider the cocycles $\omega_i$ for $1\leq i\leq 2m$.  
At least one of the following values is non-zero (depending on $\mu$):
$$\begin{sis}
&({\rm d}\omega_i)_{x_i^2}(x_ix_n,x^{\sigma-2\epsilon_i-(p-1)\epsilon_{i'}}x_n^{\mu})=
(\omega_i)_{x_i^2}([x_ix_n,x^{\sigma-2\epsilon_i-(p-1)\epsilon_{i'}}x_n^{\mu}])=\\
&\hspace{1,5cm}=(-3-\mu)(\omega_i)_{x_i^2}(x^{\sigma-\epsilon_i-(p-1)\epsilon_{i'}}x_n^{\mu})=
2(-3-\mu),\\
&({\rm d}\omega_i)_{x_i^2}(x_i^3,x^{\sigma-4\epsilon_i-(p-1)\epsilon_{i'}}x_n^{\mu})=
(\omega_i)_{x_i^2}([x_i^3,x^{\sigma-4\epsilon_i-(p-1)\epsilon_{i'}}x_n^{\mu}])=\\
&\hspace{1,5cm}=-\mu(\omega_i)_{x_i^2}(x^{\sigma-\epsilon_i-(p-1)\epsilon_{i'}}x_n^{\mu})=-2\mu.
\end{sis}$$
On the other hand, for every $g\in H^2\left(\frac{K(n)_{\geq 1}}{K(n)_{\geq d}},F\right)$,
 it holds that 
$$\begin{sis}
&(x_i^2\circ g)(x_ix_n,x^{\sigma-2\epsilon_i-(p-1)\epsilon_{i'}}x_n^{\mu})=0,\\
&(x_i^2\circ g)(x_i^3,x^{\sigma-4\epsilon_i-(p-1)\epsilon_{i'}}x_n^{\mu})=0.
\end{sis}$$
since $[x_i^2,x_ix_n]=[x_i^2,x^{\sigma-2\epsilon_i-(p-1)\epsilon_{i'}}x_n^{\mu}]=
[x_i^2,x_i^3]=[x_i^2,x^{\sigma-4\epsilon_i-(p-1)\epsilon_{i'}}x_n^{\mu}]=0.$
\end{proof}

\section{Hamiltonian algebra}

\subsection{Definition and Basic properties}

Throughout this section we fix a field $F$ of characteristic $p\neq 2,3$ and 
an even integer $n=2m\geq 2$.

We are going to use all the notations about multindices introduced at the beginning
of section 2.1. 
We are going to use often the following special $n$-tuples: $\0:=(0,\cdots,0)$,
$\sigma:=(p-1,\cdots,p-1)$ and $\sigma^i:=\sigma-(p-1)\ep_i-(p-1)\ep_{i'}$.

The vector space $A(n)=F[x_1,\cdots,x_n]/(x_1^p,\cdots,x_n^p)$, endowed with the grading
defined by ${\rm deg}(x^a)=|a|-2$, becomes a graded Lie algebra by mean of
\begin{equation*}
[x^a,x^b]=D_H(x^a)(x^b),
\end{equation*}
where  $D_H: A(n)\to W(n)={\rm Der}_F A(n)$ is defined by
\begin{equation*}D_H(f)=\sum_{j=1}^{2m} \sigma(j)D_j(f)D_{j'}=\sum_{i=1}^m\left[D_i(f)D_{i+m}-
D_{i+m}(f)D_i\right].
\end{equation*}
We denote with $H'(n)$ the quotient of $A(n)$ by the central element $1=x^{\0}$ so that 
there is an exact sequence of $H'(n)$-modules
\begin{equation}\label{H'-A-sequence}
0\to  \langle 1\rangle_F \to A(n)\to H'(n)\to 0,
\end{equation}
where $\langle 1 \rangle_F \cong F$ is the trivial $H'(n)$-module.

\begin{defi}
The Hamiltonian algebra is the derived subalgebra of $H'(n)$:
$$H(n):=H'(n)^{(1)}=[H'(n),H'(n)].$$
\end{defi}
\noindent There is an exact sequence of $H(n)$-modules (see 
\cite[Chap. 4, Prop. 4.4]{FS}):
\begin{equation}\label{H-H'-sequence}
0\to H(n) \to H'(n)\to  \langle x^{\sigma}\rangle_F \to 0,
\end{equation}
where $\langle x^{\sigma}\rangle_F\cong F$ is the trivial $H(n)$-module.

Note that the unique term of negative degree is 
$H(n)_{-1}=\oplus_{i=1}^n  \langle x_i \rangle_F$ where $x_i$ acts, via the adjoint 
action, as $D_H(x_i)=\sigma(i)D_{i'}$. The term of degree $0$ is 
$H(n)_0=\oplus_{1\leq i, j \leq n}  \langle x_ix_j \rangle_F$ and its
adjoint action on $H(n)_{-1}$ induces an isomorphism $H(2m)_0\cong \mathfrak{sp}(2m, F)$.
 
The algebra $H(n)$ admits a root space decomposition with respect to a canonical
Cartan subalgebra.

\begin{pro}\label{Cartan-H}
\begin{itemize}
\item[(a)] $T_H:=\oplus_{i=1}^m \langle x_ix_{i'}\rangle_F$
is a maximal torus of $H(n)$ (called the canonical maximal torus).
\item[(b)] The centralizer of $T_H$ inside $H(n)$ is the subalgebra
$C_H=\{x^a\: | \: a_{i'}=a_i\}$, which is hence a Cartan subalgebra (called
the canonical Cartan subalgebra). The dimension of $C_H$ is $p^m-2$.
\item[(c)] Let $\Phi_H:={\rm Hom}_{\mathbb{F}_p}(\oplus_{i=1}^n \langle x_ix_{i'}
\rangle_{\mathbb{F}_p}, \mathbb{F}_p)$, where $\mathbb{F}_p$ is the prime
field of $F$. We have a Cartan decomposition $H(n)=C_H\oplus_{\phi\in \Phi_H-\0}H(n)_{\phi}$,
where $H(n)_{\phi}=\{x^a\: | \: a_{i+m}-a_i\equiv \phi (x_ix_{i'}) \mod p\}$. 
The dimension of every $H(n)_{\phi}$, with $\phi\in \Phi_H-\0$, is $p^m$.
\end{itemize}
\end{pro}
\begin{proof}
See \cite[Chap. 4, Theo. 4.5 and 4.6]{FS}.
\end{proof}

\subsection{Proof of the Main Theorem \ref{H-finaltheorem}}

In this section, assuming the results of the next two sections, 
we give a proof of the Main Theorem \ref{H-finaltheorem}.
As a first step towards the proof, we compute the cohomology group of the second cohomology 
group of the $H(n)$-module $H'(n)$.

\begin{pro}\label{H-coho-H'}
The second cohomology group of $H'(n)$ is given by 
$$H^2(H(n),H'(n))=\bigoplus_{i=1}^n \langle {\rm Sq}(x_i)\rangle_F \bigoplus_{i<j} 
\langle \Pi_{ij} \rangle_F \bigoplus \langle \Phi\rangle_F,$$
where $\Pi_{ij}$ and $\Phi$ are the cocycles appearing in Theorem \ref{H-finaltheorem}.
\end{pro}
\begin{proof}
From the exact sequence (\ref{H'-A-sequence}) and using Propositions \ref{H-triv-2}
and \ref{H-H^2}, we get the exact sequence 
$$ 0 \to \bigoplus_{i=1}^n \langle {\rm Sq}(x_i)\rangle_F \bigoplus_{i<j} 
\langle \Pi_{ij} \rangle_F \bigoplus \langle \Phi \rangle_F 
\to   H^2(H(n),H'(n)) \stackrel{\partial}{\longrightarrow} $$
$$\to H^3(H(n),\langle 1\rangle_F) \to  H^3(H(n),A(n)).$$
We have to verify that the coboundary map $\partial$ is equal to zero, or in other words
that the cocycles which generate $H^3(H(n),\langle 1 \rangle_F)$ (see Proposition \ref{H-triv-3})
do not become zero in the group $H^3(H(n),A(n))$.

The cocycle $\Gamma_{ij}$ (for certain $i<j$, $j\neq i'$)
cannot be the coboundary of an element $h\in C^2(H(n),A(n))$.
Indeed we have that $\Gamma_{ij}(x_i^2,x_j^2,$ $x^{\sigma-(p-1)(\ep_{i'}+\ep_{j'})-
\ep_i-\ep_j})=4$ while the element 
$\d h(x_i^2,$ $x_j^2,x^{\sigma-(p-1)(\ep_{i'}+\ep_{j'})-\ep_i-\ep_j})$ cannot contain the 
monomial $1$ since the bracket of any two of the above elements vanish and all the 
three elements have degree greater or equal to $0$.

Assume now that $n\equiv -4 \mod p$ and suppose, by absurd,  
that the cocycle $\Xi$ is the coboundary 
of a cochain $f\in C^2(H(n),A(n))$. For a multindex $\0\leq a\leq \sigma$, call
$\phi(x^a)=-\phi(x^{\sigma-a})$ the coefficient of $1$ in the element $f(x^a,x^{\sigma-a})$.
Consider a triple  of elements 
$(x^a,x^b,x^c)\in H(n)\times H(n)\times H(n)$ such that $a+b+c=\sigma+\ep_k+\ep_{k'}$ 
(for a certain $k$) and ${\rm deg}(x^a)={\rm deg}(x^b)={\rm deg}(x^c)\geq 0$.
By taking the coefficient of $1$ in the equality $\Xi(x^a,x^b,x^c)=\d f(x^a,x^b,x^c)$ 
and using the relations $c_k\equiv -a_k-b_k \mod p$ and $c_{k'}=-a_{k'}-b_{k'} \mod p$,
 we get that 
$$\phi(x^a)+\phi(x^b)=\phi(x^{a+b-\ep_k-\ep_{k'}})+1 \hspace{0,5cm} \text{ if } \hspace{0,2cm} 
a_kb_{k'}-a_{k'}b_k\not\equiv 0 \mod p.
$$
By considering triples as above with ${\rm deg}(x^a)={\rm deg}(x^b)=0$, we get
the relations $2\phi(x_ix_j)=\phi(x_i^2)+\phi(x_j^2)$ and $2=\phi(x_i^2)+
\phi(x_{i'}^2)$, from which we deduce that the restriction of $\phi$ to $H(n)_0$
is determined by the values $\phi(x_i^2)$ for $1\leq i\leq m$.
Analogously, by taking ${\rm deg}(x^a)=0$ and ${\rm deg}(x^b)=1$, one gets that the 
restriction of $\phi$ to $H(n)_1$ is determined by the value $\phi(x_1^3)$ together 
with the restriction of $\phi$ to $H(n)_0$. Finally, by taking ${\rm deg}(x^a)=1$ and 
$1\leq {\rm deg}(x^b)=d\leq n(p-1)-5$, one gets that the values of $\phi$ on $H(n)_{d+1}$
are determined by the values of $\phi$ on $H(n)_1$ and on $H(n)_d$. Therefore the values 
of $\phi$ on the elements having degree $0\leq d\leq n(p-1)-4$ is determined by the 
values $\phi(x_i^2)$ for $1\leq i\leq m$ and $\phi(x_1^3)$.  
Explicitly, for an element $x^a\in H(n)$ such that $0\leq {\rm deg}(x^a)\leq n(p-1)-4$,
one gets the following formula 
$$\phi(x^a)=\left(\sum_{i=1}^n a_i-2\right)\phi(x_1^3)+\left( -a_1-2a_{1'}-\sum_{j\neq 1,1'}
\frac{3a_j}{2}+3\right) \phi(x_1^2)+ 
$$  
$$+\sum_{k=2}^m \frac{a_k-a_{k'}}{2} \phi(x_k^2)+\sum_{h=m+1}^n a_h. 
$$
Imposing the antisymmetric relation $\phi(x^{\sigma-a})=-\phi(x^a)$, we get the relation 
$$-(n+4)\phi(x_1^3)+\frac{3(n+4)}{2}\phi(x_1^2)-\frac{n}{2}=0,$$
which is impossible by the hypothesis $n\equiv -4 \mod p$ (and $p\neq 2$). 

Finally, the cocycles belonging to $H^3(H(n),H(n)_{-1};\langle 1\rangle_F)$ are not in the 
image of the coboundary map $\partial$ of above. Indeed, consider a cohomology class of 
$H^3(H(n),\langle 1\rangle_F)$ coming from $H^2(H(n),H'(n))$ and choose a representative 
$f\in Z^3(H(n),\langle 1\rangle_F)$ 
such that $f=\partial g$ where $g\in Z^2(H(n),H'(n))$. Since $g$ takes values in 
$H'(n)=A(n)_{\geq 0}$, then the cocycle $f$ vanish on the $3$-tuples of elements having 
non-negative degree. On the other hand, if $f$ belongs to 
$Z^3(H(n),H(n)_{-1};\langle 1\rangle_F)$, 
then by definition it must vanish on the $3$-tuples of elements such that at least one has 
negative degree. Putting together these two vanishings, we deduce that $f=0$.

\end{proof}

Now, using the above Proposition, we can prove the Main Theorem \ref{H-finaltheorem}. 

\begin{proof}[Proof of Theorem \ref{H-finaltheorem}]
From the exact sequence (\ref{H-H'-sequence}) and using that 
$H^1(H(n),$ $\langle x^{\sigma}\rangle_F)$ $=0$, 
we get the exact sequence 
$$0\to H^2(H(n),H(n))\to H^2(H(n),H'(n)) \to H^2(H(n),\langle x^{\sigma}\rangle_F),$$
so that we have to check which of the cocycles of the above Proposition \ref{H-coho-H'}
go to $0$ under the projection onto $H^2(H(n),\langle x^{\sigma}\rangle_F)$.
Clearly the cocycles ${\rm Sq}(x_i)$ take values in $H(n)=[H'(n), H'(n)]$ by definition.

Consider the cocycles $\Pi_{ij}\in H^2(H(n),H'(n))$. If $j\neq i,i'$ then the projection 
of $\Pi_{ij}$ onto $H^2(H(n),\langle x^{\sigma}\rangle_F)$ is $0$. Indeed $\Pi_{ij}(x^a,x^b)\subset 
\langle x^{\sigma}\rangle_F$ if and only if  
$a+b=\sigma-(p-1)\ep_i-(p-1)\ep_j+\ep_{i'}+\ep_{j'}$ but, for these pairs of elements,
it is easily checked that $D_i(x^a)D_j(x^b)-D_j(x^a)D_i(x^b)=
(a_ib_j-a_jb_i)x^{a+b-\ep_i-\ep_j}=0$.
On the other hand, if $j=i'$, then the only non-zero values of 
$\Pi_{ii'}$ are given by 
$$\Pi_{ii'}(x_ix^a,x_{i'}x^b)=x^{a+b+(p-1)\ep_i+(p-1)\ep_{i'}} \text{ for } 
a+b\leq \sigma^i.$$ 
Therefore, if $n=2$, the cocycle $\Pi_{12}$ satisfy  
$\Pi_{12}(x_1,x_2)=x^{\sigma}$ and hence it cannot be lifted to $H^2(H(n),H(n))$.
On the other hand, for $n\geq 4$, if we define $g_i\in C^1(H(n),H'(n))$
by $g_i(x^{\sigma^i})=x^{\sigma}$, then the only non-zero values
of the coboundary $\d g_i$ (for $1\leq i\leq m$) can be 
$$\begin{sis}
& \d g_i(x_ix^a,x_{i'}x^b)=-g_i([x_ix^a,x_{i'}x^b])=-x^{\sigma} \text{ if } a+b=\sigma^i,\\
& \d g_i(x_k,x^{\sigma^i})=[x_k,g_i(x^{\sigma^i})]=-\sigma(k)x^{\sigma-\ep_{k'}}
\text{ for any } 1\leq k\leq n.\\
\end{sis}$$
Therefore $\Pi_i=\Pi_{ii'}+\d g_i $ and clearly 
$\Pi_i\in H^2(H(n),H(n))$ since it vanish on the pairs $(x_ix^a,x_{i'}x^b)$
such that $a+b=\sigma^i$.

Consider now the cocycle $\Phi$. We want to prove that its projection onto $H^2(H(n),$
$\langle x^{\sigma}\rangle_F)$ vanish. From the explicit description of $\Phi$, it follows that its 
projection onto $\langle x^{\sigma}\rangle_F$ 
is given by 
$$\Phi(x^a,x^b)=\sum_{\stackrel{|\delta|=3, \delta+\widehat{\delta}<a}{\delta+\widehat{\delta}
=a+b-\delta}} \binom{a}{\delta}\binom{b}{\widehat{\delta}}\sigma(\delta)\delta!x^{\sigma},$$ 
where the above sum is set equal to $0$ if there are no elements $\delta$ verifying the 
hypothesis. Each element $\delta$ verifying the above hypothesis contributes to the 
summation with the coefficient 
$$\sigma(\delta)\delta! \binom{a}{\delta}\binom{b}{\widehat{\delta}}=-\sigma(\delta)
\delta!\binom{a}{\delta}\binom{a-\delta}{\widehat{\delta}}=-\frac{\sigma(\delta)}
{\widehat{\delta}!}\frac{a!}{(a-\delta-\widehat{\delta})!}\: ,$$
where in the first equality we substitute $b=\sigma-a+\delta+\widehat{\delta}$
and we use the relation $\binom{\sigma-c}{d}=(-1)^{|d|}\binom{c+d}{d}$ which follows
from the congruence $k!(p-1-k)!\equiv (-1)^{k+1} \mod p$ (for $0\leq k \leq p-1$). 
Now note that if a certain $\delta$ appears in the above summation, then it also 
appears its conjugate $\widehat{\delta}$ and we have that $\delta\neq \widehat{\delta}$ because
of the oddness of the degree $|\delta|$. Using the easy relations 
$\delta!=\widehat{\delta}!$ and $\sigma(\delta)=(-1)^{|\delta|}\sigma(\widehat{\delta})=
-\sigma(\widehat{\delta})$, it follows that the contributions of $\delta$ and 
$\widehat{\delta}$ are opposite and therefore the sum vanish.

\end{proof}

\subsection{Cohomology of the trivial module}

In this section we compute the second and third cohomology group of $H(n)$ with 
coefficients in the trivial module $F$.

\begin{pro}\label{H-triv-2}
The second cohomology group of the trivial module is equal to 
$$H^2(H(n),F)=
\begin{cases}
\bigoplus_{i=1}^n \langle \Omega_i\rangle_F \bigoplus \langle \Sigma \rangle_F 
& \text{ if } n\not \equiv -4 \mod p , \\
\bigoplus_{i=1}^n \langle \Omega_i\rangle_F \bigoplus \langle \Sigma \rangle_F 
\bigoplus \langle \Delta \rangle_F & \text{ otherwise,}
\end{cases} $$
where the only non-zero values of the above cocycles are
$$\begin{sis} 
&\Omega_i(x^a,x^b)= a_i & \text{ if } a+b=\sigma+\ep_i-(p-1)\ep_{i' },\\
&\Sigma(x_k,x_{k'})=\sigma(k), & \\
&\Delta(x^a,x^b)={\rm deg}(x^a) & \text{ if } a+b=\sigma.\\
\end{sis}$$
\end{pro}
\begin{proof}
Note that the cochain $\Delta$ is antisymmetric if and only if $n\equiv -4 \mod p$, because if 
$a+b=\sigma$ then ${\rm deg}(x^a)+{\rm deg}(x^b)=n(p-1)-4\equiv -n-4 \mod p$.

The verification that the above cochains are cocycles and are independent 
in $H^2(H(n), F)$ is straightforward and is left to the 
reader. We conclude by \cite[Thm 2.4]{Far}, which gives that 
$$\dim_F H^2(H(n), F)=\begin{cases}
n+1 & \text{ if } n\not \equiv -4 \mod p , \\
n+2 & \text{ otherwise.}
\end{cases} $$
\end{proof}

\noindent In order to compute $H^3(H(n), F)$, we use 
the Hochschild-Serre spectral sequence (see \cite{HS}) relative to 
the subalgebra $H(n)_{-1}<H(n)$:
\begin{equation}\label{H-HS-triv}
E_1^{r,s}=H^s(H(n)_{-1},C^r(H(n)/H(n)_{-1},F))\Rightarrow H^{r+s}(H(n),F).
\end{equation} 
For the first line of the second page of the above spectral sequence,
we have the equality
\begin{equation}\label{rel-coho}
E_2^{r,0}= H^r(H(n),H(n)_{-1};F),
\end{equation}
where $H^*(H(n),H(n)_{-1};F)$ are the relative cohomology groups of $H(n)$ 
with respect to the subalgebra $H(n)_{-1}$ with coefficients in the trivial 
module $F$ (as defined in \cite{CE}). 

Moreover, as remarked in \cite[Sec. 2.1]{Viv}, we can restrict ourself 
to consider homogeneous cohomology with respect to the maximal torus $T_H\subset H(n)$ 
(see Prop. \ref{Cartan-H}). 

\begin{pro}\label{H-triv-3}
The third cohomology group of the trivial module is equal to
$$H^3(H(n),F)=
\begin{sis}
& H^3(H(n),H(n)_{-1};F)\bigoplus_{i<j, i\neq j'}\langle \Gamma_{ij}\rangle_F 
\text{ if } n \not\equiv -4 \mod p,\\
& H^3(H(n),H(n)_{-1};F)\bigoplus_{i<j, i\neq j'}\langle \Gamma_{ij}\rangle_F 
\bigoplus \langle \Xi\rangle_F  \text{ otherwise.}
\end{sis}$$
where, by definition, the only non-zero values of the above cocycles are (for $j\neq i, i'$)  
$$\begin{sis}
&\Gamma_{ij}(x^a,x^b,x^c)=
a_ib_j-a_jb_i \hspace{0,4cm}\text{ if } a+b+c=\sigma-(p-1)\ep_{i'}-(p-1)\ep_{j'}+\ep_i+\ep_j,\\ 
& \Xi(x^a,x^b,x^c)=\sigma(k)[a_kb_{k'}-a_{k'}b_k] \hspace{0,2cm}
\text{ if } a+b+c=\sigma+\ep_k+\ep_{k'} \hspace{0,8cm} \text{ for some } k.
\end{sis}$$

\end{pro}
\begin{proof}
The verification that the above cochains are cocycles is straightforward and is left to the 
reader. In order to show that they freely generate the third cohomology group, we divide the proof
into four steps according to the spectral sequence (\ref{H-HS-triv}).

\underline{STEP I} : $(E_1^{0,3})_{\0}=H^3(H(n),F)_{\0}=0$ by homogeneity.

\underline{STEP II} : $(E_{\infty}^{1,2})_{\0}=\bigoplus_{i<j, i\neq j'}
\langle \Gamma_{ij}\rangle_F. $\\
With the notations of Lemma \ref{H-triv-(1,2)} below, 
consider a cochain $\zeta=\sum_{k=1}^n d_k \zeta_k
\in (E_1^{1,2})_{\0}$ and suppose that it can be lifted to a global cocycle in $Z^3(H(n),F)_{\0}$
(which we continue to call $\zeta$). Consider the following cocycle condition
$$\d \zeta(x_i,x_{i'},x_i^2,x_{i'}^2)=-4\sigma(i)\zeta(x_i,x_{i'},x_{i}x_{i'})=
-4\sigma(i)[\sigma(i')d_i-\sigma(i)d_{i'}]=4[d_i+d_{i'}],$$
from which we deduce the relation $d_i=-d_{i'}$. 
It is easily checked that $\zeta$ is the coboundary of the cocyle 
$f\in (E_1^{0,2})_{\0}=H^2(H(n)_{-1},F)_{\0}$ defined by 
$f(x_i,x_{i'})=-d_{i'}=d_i$, since we have (for $j\neq i, i'$) 
$$\zeta(x_ix_j,x_{i'},x_{j'})=\sigma(j')d_{i'}-\sigma(i')d_{j'}=
{\rm d}(f)(x_ix_j,x_{i'},x_{j'}).$$

Suppose now that $n\geq 4$.
The cocycles $\overline{\Gamma_{ij}}$ with $j \neq i,i'$ appearing in  
Lemma \ref{H-triv-(1,2)} are clearly lifted by the cocycles $\Gamma_{ij}$.
On the other hand, the cocycles $\overline{\Gamma_{ii'}}$ cannot be lifted to $Z^3(H(n),F)_{\0}$.
Indeed, by absurd, suppose that we can find such a lift and call it $\Gamma_{ij}
\in Z^3(H(n),F)_{\0}$.
We can suppose that $\Gamma_{ij}$ takes its non-zero values on the triples 
$(x^{\alpha},x^{\beta},x^{\gamma})$ such that $\alpha+\beta+\gamma=\sigma^i+\ep_i+\ep_{i'}$,
where $\sigma^i:=\sigma-(p-1)\ep_i-(p-1)\ep_{i'}$. 
Consider the following cocycle condition (where $a, b, c$ are multindices verifying 
$a+b+c=\sigma^i$):
$$0=\sigma(i)\Gamma_{ii'}(x_i^2x^a,x_{i'},x_{i'}x^b,x^c)=-2\Gamma_{ii'}(x_ix^a,x_{i'}x^b,x^c)
-2\Gamma_{ii'}(x_ix^{a+b},x_{i'},x^c).$$
We deduce that the value of $\Gamma_{ii'}(x_ix^a,x_{i'}x^b,x^c)$ depends only on the multindex 
$c$ and therefore, for every $\0\leq c\leq \sigma^i$, we can define 
$\omega(c):=\Gamma_{ii'}(x_ix^a,x_{i'}x^b,x^c)$ for every pair of 
indices $a, b$ such that $a+b+c=\sigma^i$. By the fact  that 
${\Gamma_{ij}}_{|H(n)_{-1}\times H(n)_{-1}}=\overline{\Gamma_{ij}}$, we get 
$$\begin{sis}
&\omega(\sigma^i)=\Gamma_{ij}(x_i,x_{i'},x^{\sigma^i})=\overline{\Gamma_{ij}}
(x_i,x_{i'},x^{\sigma^i})=1, \\
& \omega(\epsilon_j)=\Gamma_{ij}(x_ix^{\sigma^i-\ep_j},x_{i'},x_j)=
\overline{\Gamma_{ij}}(x_ix^{\sigma^i-\ep_j},x_{i'},x_j)=0.
\end{sis}$$
Finally consider the following cocycle condition where $j\neq i,i'$ and $\0\leq d\leq \sigma^i$
is a multindex such that $d_{j'}>0$:
$$0=\sigma(j)\Gamma_{ii'}(x_i,x^{\sigma^i-d+\ep_{j'}}x_{i'},x_j,x^d)=
d_{j'}[\omega(d)-\omega(d-\ep_{j'})],
$$
where we used that $\omega(\epsilon_j)=0$. We deduce that the value  $\omega(d)$ does not depend 
on the coefficient $d_{j'}$ and, by repeating for every index $j\neq i,i'$, we conclude
that $\omega$ must be constant. But this contradicts with $\omega(\epsilon_j)=0$ and 
$\omega(\sigma^i)=1$.

\underline{STEP III} : $(E_{\infty}^{2,1})_{\0}=
\begin{cases}
\langle \Xi \rangle_F & \text{ if } n\equiv -4 \mod p,\\
0 & \text{ otherwise.}
\end{cases}$ \\ 
With the notations of Lemma \ref{H-triv-(2,1)} below, consider the cochain
 $\xi=\sum_{k=1}^n e_k\xi_k$ and suppose that it can be lifted to a global cocycle of
$Z^3(H(n),F)_{\0}$ (which, as usual, we continue to call $\xi$). 
From the following cocycle condition
$$0=\d \xi(x_{k'},x_k^2,x_{k'}^2,x^{\sigma-\ep_{k'}})=\sigma(k')
\xi(x_k^2,x_{k'}^2,x^{\sigma-\ep_k-\ep_{k'}})+4e_k$$
together with the analogues one obtained interchanging $k$ with $k'$, we get that
$e_k=e_{k'}$. 

If $n \equiv -4 \mod p$, then the cochain $\sum_{i=1}^n \xi_i$ 
is lifted by the global cocycle $\Xi\in (E_{\infty}^{2,1})_{\0}$.
We will show that under the assumption that either 
$n\not\equiv -4 \mod p$ or $n\equiv -4 \mod p$ and $\sum_{i=1}^m e_i=0$, then
$\xi$ belongs to the image of the differential map 
$$\d:(E_1^{1,1})_{\0}\to (E_1^{2,1})_{\0},
$$
coming from the spectral sequence (\ref{H-HS-triv}).

Consider the cochains $\eta_k \in (E_1^{1,1})_{\0}=
H^1(H(n)_{-1},C^1(H(n)/H(n)_{-1},F))_{\0}$ (for $1\leq k\leq m$), whose only non-zero
values are given by 
$$ \eta_k(x_k, x^{\sigma-\ep_k})=\eta_k(x_{k'},x^{\sigma-\ep_{k'}})=-1. $$
Form the cochain $\eta:=\sum_{i=1}^{m} (\frac{e_i}{2}+\beta)\eta_i$, where 
$\beta \in F$ is defined as
$$\beta=\begin{sis}
-\frac{\sum_{i=1}^m e_i}{n+4} & \text{ if } n\not\equiv -4 \mod p\\
0 & \text{ if } n\equiv -4 \mod p \text{ and }\sum_{i=1}^m e_i=0.
\end{sis}$$
It is straigthforward to check that the cocycle 
$\xi-\d \eta \in C^1(H(n)_{-1},C^2(H(n)/$ $H(n)_{-1},F))_{\0}$
is the coboundary of the cochain $g\in C^2(H(n)$ $/H(n)_{-1},F)$
defined by (and vanishing elsewhere)
$$g(x^a,x^{\sigma-a})= \sum_{i=1}^m (a_i+a_{i'})\frac{e_i}{2}+\deg_H(x^a)\beta
\quad \text{ if } |a|, |\sigma-a| \geq 2.  $$
This shows that $[\xi]=[\d \eta]\in (E_1^{2,1})_{\0}$.


Suppose next that the cocycle $\rho_{ij}$ (for certain $i<j$) can be lifted to a 
global cocycle of $Z^3(H(n),F)$, which we continue to call $\rho_{ij}$. 
For $l=2,\cdots,p-1$, we define $f_l:=\rho_{ij}(x_ix_{j'},x^{\sigma-l\ep_{j'}},
x^{\sigma-(p-1)\ep_{i'}-(p+1-l)\ep_{j'}})$.
Consider the following cocycle condition for $1\leq l\leq p-1$:
$$0=\d \rho_{ij}(x_ix_{j'},x_j,x^{\sigma-l\ep_{j'}},x^{\sigma-(p-1)\ep_{i'}-
(p-l)\ep_{j'}})=$$
\begin{equation*}
=\frac{-(1+\delta_{ji'})}{l}+\delta_{1l}(1+\delta_{ji'})-\sigma(j)(l+1)f_{l+1}+
\sigma(j)(l-1)f_l.\tag{*}
\end{equation*}
The above equation (*) with $l=1,\cdots,p-2$ gives that 
$$f_l=\frac{-(1+\delta_{ji'})(l-2)\sigma(j)}{(l-1)l} \hspace{1cm}
\text{ for } l=2,\cdots,p-1. $$
Substituting in the above equation (*) with $l=p-1$, we get
$$0=\frac{-(1+\delta_{ji'})}{p-1}+\sigma(j)(p-2)\frac{-(1+\delta_{ji'})(p-3)\sigma(j)}
{(p-2)(p-1)}=(1+\delta_{ji'})-3(1+\delta_{ji'}),$$
which is impossible since $p\neq 2$.
Therefore the cocycles $\rho_{ij}$ do not belong to $(E_{\infty}^{2,1})_{\0}$.

\underline{STEP IV} : $(E_{\infty}^{3,0})_{\0}=(E_2^{3,0})_{\0}=H^3(H(n),H(n)_{-1};F)_{\0}$.

From Proposition \ref{H-triv-2}, one can easily deduce that 
$$\begin{sis}
&(E_{\infty}^{0,2})_{\0}=(E_2^{0,2})_{\0}=\langle \Sigma \rangle_F,\\
&(E_{\infty}^{1,1})_{\0}=(E_2^{1,1})_{\0}=
\begin{cases}
\bigoplus_{i=1}^n \langle \Omega_{i}\rangle_F \oplus \langle \Delta \rangle_F &
\text{ if } n \equiv -4 \mod p, \\   
\bigoplus_{i=1}^n \langle \Omega_{i}\rangle_F & \text{ otherwise. }
\end{cases}\\
\end{sis}$$
This implies that $(E_{\infty}^{3,0})_{\0}=(E_2^{3,0})_{\0}$ and the result follows
from equality (\ref{rel-coho}).

\end{proof}

\begin{rem}
It can be proved that $H^3(H(n),H(n)_{-1};F)=\oplus_{i=1}^n \langle \Upsilon_i
\rangle_F$ where the cocycles $\Upsilon_i$ are defined by 
$$\Upsilon_i(x^a,x^b,x^c)=\sigma(k)[a_kb_{k'}-a_{k'}b_k] \text{ if } a+b+c=\sigma+
p\ep_i+\ep_k+\ep_{k'} \text{ for some } k.$$
We omit the proof, since we do not need this result to prove the Main Theorem
\ref{H-finaltheorem}.
\end{rem}

\begin{lem}\label{H-triv-(1,2)}
In the above spectral sequence (\ref{H-HS-triv}), we have that
$$(E_1^{1,2})_{\0}= 
\begin{sis}
&\frac{\bigoplus_{k=1}^{n}  \langle \zeta_k \rangle_F}{\langle \sum_{k=1}^n 
\sigma(k)\zeta_k\rangle_F}
\bigoplus_{i<j}  \langle \overline{\Gamma_{ij}} \rangle_F & \text{ if } n\geq 4,\\
& \frac{\bigoplus_{k=1}^n \langle \zeta_k \rangle_F}{\langle \sum_{k=1}^n \sigma(k)
\zeta_k\rangle_F}
& \text{ if } n=2, 
\end{sis}$$
where the only non-zero values of the above cocycles are
$$\begin{sis}
&\overline{\Gamma_{ij}}(x_i,x_j,x^{\sigma-(p-1)\ep_{i'}-(p-1)\ep_{j'}})=1,\\
& \zeta_k(x_k,x_h,x_{k'}x_{h'})=\sigma(h) \: \text{ for } \: h=1,\cdots, n.
\end{sis}$$
\end{lem}
\begin{proof}
Consider the exact sequence 
\begin{equation}\label{H-exa1}
0\to C^1(x^{\sigma},F) \to C^1(H'(n)/H'(n)_{-1},F) \to C^1(H(n)/H(n)_{-1},F)\to 0,
\end{equation}
where $C^1(x^{\sigma},F)$ is a trivial $H(n)_{-1}$-module. The coboundary map
$$ H^1(H(n)_{-1},C^1(H(n)/H(n)_{-1},F))_{\0}
\stackrel{\partial^{(2)}}{\longrightarrow} H^2(H(n)_{-1},C^1(x^{\sigma},F))_{\0}$$
is surjective. Indeed consider the cocycles 
$\eta_k\in H^1(H(n)_{-1},C^1(H(n)/H(n)_{-1},F))_{\0}$ (with $1\leq k\leq m$), defined as 
$$\eta_k(x_k, x^{\sigma-\ep_k})=\eta_k(x_{k'},x^{\sigma-\ep_{k'}})=-1.$$
It is easy to check that $\partial^{(2)}$ sends $\eta_k$ into the cocycles 
$\{ (x_k,x_{k'},x^{\sigma})\mapsto -2 \}$ which generate the last group 
$H^2(H(n)_{-1},C^1(x^{\sigma},F))_{\0}$.

Using the above surjectivity, together with the vanishing 
$H^3(H(n)_{-1},C^1(x^{\sigma},$ $ F))_{\0}=0$ which follows directly by homogeneity
considerations, we get that
$$(E_2^{1,2})_{\0}=H^2(H(n)_{-1},C^1(H'(n)/H'(n)_{-1},F))_{\0}.$$
Consider now the following exact sequence of $H(n)_{-1}$-modules 
\begin{equation}\label{H-exa2}
0 \to C^1(H'(n)/H'(n)_{-1},F) \to C^1(A(n),F) \to C^1(A(n)_{<0},F) \to 0,
\end{equation}
obtained from the fact that $H'(n)/H'(n)_{-1}=A(n)/A(n)_{<0}$ (see (\ref{H'-A-sequence})).
Using the following isomorphism of $H(n)_{-1}$-modules 
\begin{equation}\label{H-isomo}
\begin{aligned}
\chi: A(n) & \stackrel{\cong}{\longrightarrow} C^1(A(n),F)\\
      x^a & \mapsto \chi_{x^a}(x^b)=
\begin{sis}
1 & \text{ if } b=\sigma-a,\\
0 & \text{ otherwise, } 
\end{sis}
\end{aligned}
\end{equation}
together with \cite[Prop. 3.4]{Viv}, we get that 
$$H^1(H(n)_{-1},C^1(H'(n)/H'(n)_{-1},F))_{\0}= H^1(H(n)_{-1},C^1(A(n),F))_{\0}=
\oplus_{i=1}^n  \langle \overline{\Omega_i} \rangle_F,$$ where the cocycles 
$\overline{\Omega_i}$ are defined by (and vanish outside)
$\overline{\Omega_i}(x_i,x^{\sigma-(p-1)\ep_{i'}})=1.$
Therefore we get the exact sequence
$$0 \to  H^1(H(n)_{-1},C^1(A(n)_{<0},F))_{\0} \stackrel{\partial^{(2)}}{\longrightarrow}
(E_1^{1,2})_{\0}\to H^2(H(n)_{-1},C^1(A(n),F))_{\0}.$$
The first group on the left is generated over $F$ by the cocycles 
$\widetilde{\zeta}_k$ ($k=1,\cdots, n$) defined by $\widetilde{\zeta_k}(x_k,x_{k'})
=1$ and subject to the relation $\sum_{k=1}^n \sigma(k')\widetilde{\zeta}_k=0$
coming from the element $\langle 1\mapsto 1 \rangle_F\in C^1(A(n)_{<0},F)_{\0}$.
It is easily checked that $\partial^{(2)}(\widetilde{\zeta}_k)=\zeta_k$. 

Moreover, using the isomorphism (\ref{H-isomo}) of $H(n)_{-1}$-modules 
$A(n)\cong C^1(A(n),F)$  and \cite[Prop. 3.4]{Viv}, we get that 
$H^2(H(n)_{-1},C^2(A(n),F))_{\0}$ is freely generated over $F$ by the cocycles 
$\overline{\Gamma_{ij}}$ for $1\leq i< j\leq n$. We conclude by observing that 
$\overline{\Gamma_{ij}}$ can be lifted to
$H^2(H(n)_{-1},C^1(H'(n)/H'(n)_{-1},F))_{\0}$ if and only if $n\geq 4$.
\end{proof}

\begin{lem}\label{H-triv-(2,1)}
In the above spectral sequence (\ref{H-HS-triv}), we have that
$$(E_1^{2,1})_{\0}=
\begin{cases}
\bigoplus_{k=1}^n \langle \xi_k \rangle_F \bigoplus_{i<j}\langle 
\rho_{ij}\rangle_F 
& \text{ if } n\geq 4, \\
\bigoplus_{k=1}^n \langle \xi_k \rangle_F \oplus  
& \text{ if } n=2, \\  
\end{cases}$$
where the only non-zero values of the above cocycles are 
$$\begin{sis}
& \xi_k(x_h,x_{h'}x_k,x^{\sigma-\ep_k})= \sigma(h)(1+\delta_{h'k}) 
&\text{ for every } h=1,\cdots, n ,\\ 
& \rho_{ij}(x_i,x^{\sigma-l\ep_{j'}},x^{\sigma-(p-1)\ep_{i'}-(p-l)\ep_{j'}})
=-\frac{\sigma(j)}{l} & \text{ for every } l=1,\cdots, p-1,\\
&\rho_{ij}(x_j,x^{\sigma-l\ep_{i'}},x^{\sigma-(p-l)\ep_{i'}-(p-1)\ep_{j'}})
=\frac{\sigma(i)}{l} & \text{ for every } l=1,\cdots, p-1. 
\end{sis}$$
\end{lem}
\begin{proof}
Consider the following exact sequence of $H(n)_{-1}$-modules
\begin{equation}\label{H-exa3}
0\to C^2(H'(n)/H'(n)_{-1},F) \to C^2(A(n),F) \stackrel{{\rm res}}{\longrightarrow}
 C^1(A(n)_{<0} \times A(n),F)\to 0,
\end{equation} 
It is easy to see
that $C^1(A(n)_{<0}\times A(n),F)_{\0}^{H(n)_{-1}}$ is generated by the cocycle $\zeta$ defined 
by $\zeta(1,x^{\sigma})=\zeta(x_i,x^{\sigma-\ep_i})=1$ (for every $i=1,\cdots, n$) and that the 
image of $\zeta$ under the first coboundary map is non-zero and equal to
$-\sum_{k=1}^n \xi_k$.
Therefore, using the Lemma \ref{H-vanish-H^1} below, we get that  
\begin{equation}\label{H-equa-tec1}
H^1(H(n)_{-1},C^2(H'(n)/H'(n)_{-1},F))_{\0}=
\langle \sum_{k=1}^n \xi_k \rangle_F.
\end{equation}
Consider finally the following exact sequence 
\begin{equation}\label{H-exa4}
C^1(H(n)/H(n)_{-1},F) \stackrel{\theta}{\hookrightarrow} C^2(H'(n)/H'(n)_{-1},F)
\twoheadrightarrow C^2(H(n)/H(n)_{-1},F) ,
\end{equation}
where the map $\theta$ sends the cocycle $g$ into the cocycle $\theta(g)$ defined 
by $\theta(g)(x^{\sigma},x^a)=g(x^a)$.
By taking cohomology, we get the exact sequence
$$ H^1(H(n)_{-1},C^2(H'(n)/H'(n)_{-1},F))_{\0}\to (E_1^{2,1})_{\0}
\stackrel{\partial^{(2)}}{\longrightarrow} (E_1^{1,2})_{\0}.$$
We conclude by using (\ref{H-equa-tec1}), Lemma \ref{H-triv-(1,2)} and the facts 
that $\partial^{(2)}(\rho_{ij})=2
\overline{\Gamma_{ij}}$ and $\partial^{(2)}(\xi_k)=\sigma(k)\zeta_{k'}$.

\end{proof}

\begin{lem}\label{H-vanish-H^1}
Consider $A(n)$ as a $H(n)_{-1}$-module. Then we have that
$$H^1(H(n)_{-1},C^2(A(n),F))_{\0}=0.$$
\end{lem}
\begin{proof}
During this proof, we use the generators $D_i:=\sigma(i') x_{i'}$ of $H(n)_{-1}$. 
Moreover, if $g\in C^2(A(n),1)$, we set $\widetilde{g}(x^a,x^b):=\frac{g(x^a,x^b)}{a!b!}$ 
where, as usual, for a multindex $a=(a_1,\cdots, a_n)$ we set $a!:=\prod_i a_i!$.
Analogously, if $f\in C^1(H(n)_{-1},C(A(n),F))$, we set 
$\widetilde{f_{D_i}}(x^a,x^b):=\frac{f_{D_i}(x^a,x^b)}{a!b!}$ where $f_{D_i}\in 
C^2(A(n),F)$ denotes, as usual, the value of $f$ on $D_i\in H(n)_{-1}$.

Take a homogeneous cochain $f\in Z^1(H(n)_{-1},C^2(A(n),F))_{\0}$
The cocycle conditions
for $f$ are $D_i\circ f_{D_j}=D_j\circ f_{D_i}$ for every $1\leq i, j\leq n$. 

\underline{STEP I}: The cocycle $f$ verifies the following condition 
\begin{equation*}
0=\sum_{k=-\alpha_i}^{\beta_i} (-1)^k \widetilde{f_{D_i}}(x^{\alpha+k\ep_i},x^{\beta-k\ep_i})
\text{ for every } (x^{\alpha},x^{\beta}) \text{ such that } 
\alpha_i+\beta_i\leq p-1 .
\tag{*}\end{equation*}
For every pair $(x^{\alpha},x^{\beta})$ as above (such that 
$\alpha_i+\beta_i\leq p-1$), define
$$\phi_i(x^{\alpha},x^{\beta}):=\sum_{k=-\alpha_i}^{\beta_i} 
(-1)^k \widetilde{f_{D_i}}(x^{\alpha+k\ep_i},x^{\beta-k\ep_i}).$$
We have to prove that $\phi_i(x^{\alpha},x^{\beta})=0$. Using the cocycle conditions
$D_j\circ \widetilde{f_{D_i}}=D_i\circ \widetilde{f_{D_j}}$ and a telescopic sum, 
it is easy to see that $(D_j\circ \phi_i)(x^{\alpha},x^{\beta})=0$ for every $j\neq i$.
From these conditions, we get that (for every index
$j\neq i$)
$$\phi_i(x^{\alpha},x^{\beta})=
\begin{sis}
&0  \hspace{4cm} \text{ if } \alpha_j+\beta_j<p-1 \text{ for some } j\neq i,\\
&(-1)^{|\beta|-\beta_i}\phi_i(x^{\alpha+\beta-\sigma+(p-1-\beta_i)\ep_i},
x^{\sigma-(p-1-\beta_i)\ep_i}) \hspace{0,5cm} \text{ otherwise. }\\
\end{sis}$$
So assume we are in the second case, that is $\alpha_j+\beta_j\geq p-1$ for every
$j\neq i$. Consider the same formula of above for the couple $(x^{\beta},x^{\alpha})$.
 By using using the antisymmetry of $\phi_i$ and the property $\phi_i(x^{\alpha+d\ep_i},
x^{\beta-d\ep_i})=(-1)^d \phi_i(x^{\alpha},x^{\beta})$ for $-\alpha_i\leq d\leq p-1-\alpha_i$ 
and $-\beta_i\leq d\leq p-1-\beta_i$, we obtain
$$\left[(-1)^{|\beta|}+(-1)^{|\alpha|}\right]
\phi_i(x^{\alpha+\beta-\sigma+(p-1-\beta_i)\ep_i},x^{\sigma-(p-1-\beta_i)\ep_i})=0.$$
Now recall that $f$ is homogeneous and therefore we have to consider only the 
pairs $(x^{\alpha},x^{\beta})$ such that the sum of the weights of $x^{\alpha}$, 
$x^{\beta}$ and $D_i$ is $\0$. Using the conditions
$\alpha_i+\beta_i\leq p-1$ and $\alpha_j+\beta_j\geq p-1$ for every 
$j\neq i$, we find the equalities (and not merely the congruences modulo $p$):
$$\alpha_{i'}+\beta_{i'}=\alpha_i+\beta_i+(p-1) \: \text{ and } \:\alpha_{j'}+\beta_{j'}=
\alpha_j+\beta_j \text{ for every } j\neq i,i'.$$
We deduce that $|\alpha|+|\beta|$ is even and, substituting in the expression above,
we get the required vanishing.

\underline{STEP II}: The cocycle $f$ is a coboundary.\\
We have to find an element $g\in C(A(n),F)_{\0}$ such that $f_{D_i}=D_i\circ g$. 
For a homogeneous pair $(x^a,x^b)$ (that is a pair such the sum of the weights of $x^a$ 
and $x^b$ is $\0$), we define
\begin{equation*}
\widetilde{g}(x^a,x^b)=
\begin{sis}
&\hspace{0,6cm} \sum_{k=0}^{a_i} (-1)^k \widetilde{f_{D_i}}(x^{a-k\ep_i},x^{b+(k+1)\ep_i})
\hspace{1,4cm} \text{ if } a_i+b_i<p-1, \\
&\sum_{\emptyset\neq I\subset\{1,\cdots,n\}} 
\sum_{(c,d)\in S_I(a,b)}
\frac{(-1)^{|I|-1}{\rm sign}(c,d)}{2}
\widetilde{f_{D_I}}(x^c,x^d)  \text{ if } a+b\geq \sigma,
\end{sis}
\end{equation*}
where for a non-empty subset $I$ of $\{1,\cdots, n\}$ (of cardinality $|I|$), 
we define $S_I(a,b)$ to be the set of pairs $(c,d)$ of multindices verifying: $c_i+d_i=
a_i+b_i+1$ and ${\rm min}(a_i,b_i)+1\leq c_i,d_i \leq {\rm max}(a_i,b_i)$ if $i\in I$
and $c_j=a_j$, $d_j=b_j$ if $j\not\in I$ (in particular $S_I(a,b)\neq \emptyset$
if and only if $a_i\neq b_i$ for $i\in I$). Moreover, if $(c,d)\in S_I(a,b)$, we put 
${\rm sign}(c,d)=\prod_{i\in I}{\rm sign}_i(c,d)$ and ${\rm sign}_i(c,d)=d_i-b_i$ or 
$c_i-a_i$ according, respectively, to the cases $b_i<a_i$ and $a_i<b_i$.
Finally, if $I=\{i_1,\cdots,i_r\}$, we put $f_{D_I}:=D_{i_1}\circ \cdots\circ D_{i_{r-1}}\circ
f_{D_{i_r}}$ which does not depend upon the order of the elements of $I$ by the 
cocycle conditions verified by $f$.

The cochain $g$ is well-defined because if $i$ and $j$ are two indices such that $a_i+b_i<p-1$
and $a_j+b_j<p-1$, then the following expression 
$$\widetilde{g}(x^a,x^b)=\sum_{k=0}^{a_i}\sum_{h=0}^{a_j}(-1)^{k+h}(D_i\circ 
\widetilde{f_{D_j}})(x^{a-k\ep_i-h\ep_j},x^{b+(k+1)\ep_i+(h+1)\ep_j})$$
is symmetric in $i$ and $j$ because of $D_i\circ f_{D_j}=D_j\circ f_{D_i}$
and reduces, via a telescopic sum, to the first expression occuring in the definition of 
$\widetilde{g}$.

Moreover it is clear from the definition that $\widetilde{g}$ is antisymmetric 
in the case $a+b\geq \sigma$, while in the case $a_i+b_i<p-1$ (for a certain $i$) the 
antisymmetry follows from the condition $(*)$ of above. 

Finally we have to check that $(D_i\circ \widetilde{g})(x^{\alpha},x^{\beta})=
\widetilde{f_{D_i}}(x^{\alpha},x^{\beta})$ for every index $i$ and 
every pair $(x^{\alpha},x^{\beta})$ such that 
the sum of the weights of $x^{\alpha}$, $x^{\beta}$ and of $D_i$ is $\0$. If $\alpha_i=
\beta_i=0$ then $(D_i\circ \widetilde{g})(x^{\alpha},x^{\beta})=0$ and 
$\widetilde{f_{D_i}}(x^{\alpha},x^{\beta})=0$ by the 
condition (*) of above. If $\alpha_i=0$ and $\beta_i<0$ then we get that 
$(D_i\circ \widetilde{g})(x^{\alpha},x^{\beta})=-\widetilde{g}(x^{\alpha},x^{\beta-\ep_i})$ 
is equal to $\widetilde{f_{D_i}}(x^{\alpha},x^{\beta})$ by the first case of the 
definition of $\widetilde{g}$.
The case $\alpha_i>0$ and $\beta_i=0$ follows from the preceding one by the antisymmetry of 
$g$. Therefore we are left  with the case $\alpha_i,\beta_i>0$. 

Suppose first that $\alpha+\beta-\ep_i\not \geq \sigma$. Take an index $j$ (may be equal to
$i$) such that $(\alpha+\beta-\ep_i)_j<p-1$. Using the first case of the definition of $g$,
we have
$$(D_i\circ \widetilde{g})(x^{\alpha},x^{\beta})=\sum_{k=0}^{\alpha_j}(-1)^{k+1}
(D_i\circ \widetilde{f_{D_j}})(x^{\alpha-k\ep_j},x^{\beta+(k+1)\ep_j})=$$
$$=\sum_{k=0}^{\alpha_j}(-1)^{k+1}
(D_j\circ \widetilde{f_{D_i}})(x^{\alpha-k\ep_j},x^{\beta+(k+1)\ep_j})=
\widetilde{f_{D_i}}(x^{\alpha},x^{\beta}),
$$
where in the last equality we used a telescopic summation.

On the other hand, suppose that $\alpha+\beta-\ep_i\geq \sigma$. We need two auxiliary 
facts before proving the required equality in this case. First of all, observe  
that the hypothesis $\alpha+\beta-\ep_i\geq \sigma$
forces the equalities (and not merely the congruences modulo $p$) $\alpha_i+\beta_i-1=
\alpha_{i'}+\beta_{i'}$ and $\alpha_j+\beta_j=\alpha_{j'}+\beta_{j'}$ for every 
$j\neq i,i'$. Therefore the sum of the degrees of the 
multindices $|\alpha|+|\beta|$ must be odd. 
Moreover, we can re-write the second expression occurring in the definition of
$g$ in a way that will be more suitable for our purpose. 
Fix an index $i$, a homogeneous pair $(x^a,x^b)$ satisfying $a+b\geq \sigma$ 
and suppose that $a_i<b_i$. By splitting the summation occurring in the definition
of $\widetilde{g}(x^a,x^b)$ according to the cases $I=\{i\}$, $I=\{i\}\cup J$ and $I=J$
with $i\in J\neq \emptyset$, and using a telescopic summation, we get  
\begin{equation*}
2\widetilde{g}(x^a,x^b)=\sum_{k=1}^{b_i-a_i}(-1)^k\widetilde{f_{D_i}}(x^{a+k\ep_i},
x^{b+(1-k)\ep_i})+ \tag{**}
\end{equation*}
\begin{equation*}
+\sum_{i\not\in J\neq \emptyset}\sum_{(c,d)\in S_J(a,b)}
(-1)^{|J|+b_i-a_i+1}{\rm sign}(c,d)\widetilde{f_{D_J}}(x^{c+(b_i-a_i)\ep_i},
x^{d-(b_i-a_i)\ep_i}). 
\end{equation*}
If $a_i=b_i$ then the above expression is trivially true while if $a_i>b_i$ 
then we get an analogous expression using the antisymmetry of $g$.

Finally, in order to prove the required equality $(D_i\circ \widetilde{g})(x^{\alpha},x^{\beta})=
\widetilde{f_{D_i}}(x^{\alpha},x^{\beta})$, we have to distinguish two cases: 
$\alpha_i<\beta_i-1$ and $\alpha_i=\beta_i$ (the case $\alpha_i>\beta_i$ follows by antisymmetry).
In the first case $\alpha_i<\beta_i$, consider
$$-2(D_i\circ \widetilde{g})(x^{\alpha},x^{\beta})=
2\widetilde{g}(x^{\alpha-\ep_i},x^{\beta})+2\widetilde{g}(x^{\alpha},x^{\beta-\ep_i})
$$
and apply formula $(**)$ to the terms $( x^{\alpha-\ep_i},x^{\beta})$ and 
$(x^{\alpha},x^{\beta-\ep_i})$, which verify the required conditions in virtue 
of our hypothesis. By summing the first terms in the corresponding expressions 
$(**)$, we get
\begin{equation*}
-\widetilde{f_{D_i}}(x^{\alpha},x^{\beta})+(-1)^{\delta_i}
\widetilde{f_{D_i}}(x^{\alpha+\delta_i\ep_i},
x^{\beta-\delta_i\ep_i}),\tag{***1}
\end{equation*} 
where we put $\delta_i:=\beta_i-\alpha_i-1\geq 0$.
By summing the last terms in the corresponding expressions $(**)$ and using that 
if $i\not\in J$ then a pair $(c,d)$ belongs to $S_J(\alpha,\beta)$ if and only if 
$(c-\ep_i,d)\in S_J(\alpha-\ep_i,\beta)$ (and, analogously, if and only if 
$(c, d-\ep_i)\in S_J(\alpha,\beta-\ep_i)$), we obtain
\begin{equation*}
\sum_{i\not\in J\neq \emptyset} \sum_{(c,d)\in S_J(\alpha,\beta)}
(-1)^{\delta_i+|J|}{\rm sign}(c,d)(D_i\circ \widetilde{f_{D_J}})
(x^{c+\delta_i\ep_i},x^{d-\delta_i\ep_i}).\tag{***2}
\end{equation*}
Using that $D_i\circ \widetilde{f_{D_J}}=D_J\circ \widetilde{f_{D_i}}$ 
and iterated telescopic summations, the above expression $(***2)$
reduces to 
\begin{equation*}
-(-1)^{\delta_i}\widetilde{f_{D_i}}(x^{\alpha+\delta_i\ep_i},x^{\beta-\delta_i\ep_i})
-(-1)^{|\beta|-|\alpha|}\widetilde{f_{D_i}}(x^{\beta},x^{\alpha}). \tag{***2'}
\end{equation*}
Summing expressions $(***1)$ and $(***2')$ and using the fact that $|\beta|+|\alpha|$
is odd, we end up with $-\widetilde{f_{D_i}}(x^{\alpha},x^{\beta})+
\widetilde{f_{D_i}}(x^{\beta},x^{\alpha})=
-2\widetilde{f_{D_i}}(x^{\alpha},x^{\beta})$, q.e.d.

The proof in the other case $\alpha_i=\beta_i$ is similar apart from the 
fact that one has to use both the expression $(**)$ and the analogous one with $a_i>b_i$. 
We leave the details to the reader.
\end{proof}

\subsection{Cohomology of $A(n)$}

In this section we compute the second cohomology group 
of the $H(n)$-module $A(n)$.

\begin{pro}\label{H-H^2}
The second cohomology group of $A(n)$ is given by  
$$H^2(H(n),A(n))=
\begin{sis}
& \bigoplus_{i=1}^n \langle {\rm Sq}(x_i) \rangle_F \bigoplus_{i=1}^n 
\langle \Omega_i\rangle \bigoplus_{i<j} \langle \Pi_{ij} \rangle_F \bigoplus \langle 
\Phi \rangle_F \text{ if } n\not\equiv -4\mod p,\\
&\bigoplus_{i=1}^n \langle {\rm Sq}(x_i) \rangle_F \bigoplus_{i=1}^n 
\langle \Omega_i\rangle \bigoplus_{i<j} \langle \Pi_{ij} \rangle_F \bigoplus \langle 
\Phi \rangle_F\bigoplus \langle \Delta\rangle_F \text{ otherwise. }
\end{sis}$$
where $\Omega_i$ and $\Delta$ are the cocycles of Proposition \ref{H-triv-2}, 
$\Phi$ and $\Pi_{ij}$ (with 
$j\neq i'$) are the cocycles of Theorem \ref{H-finaltheorem} and the remaining  cocycles 
$\Pi_{ii'}$ are defined by (and vanish outside): 
$$
 \Pi_{ij}(x_ix^a,x_{i'}x^b)=x^{a+b+(p-1)\ep_i+(p-1)\ep_{i'}} 
\hspace{1cm} \text{ if } a+b\leq \sigma^i.
$$ 
\end{pro}
\begin{proof}
It is straightforward to verify that the above cochains are cocyles
and that they are independent in $H^2(H(n), A(n))$.
Therefore it is enough to prove that 
$$\dim_F H^2(H(n),A(n))=
\begin{sis}
& \binom{n}{2}+2n+1 & \text{ if } n\not\equiv -4\mod p,\\
&\binom{n}{2}+2n+2 & \text{ otherwise. }
\end{sis}$$ 

\noindent It is easily seen that $A(n)$ is the restricted $H(n)$-module induced
from the restricted trivial $H(n)_{\geq 0}$-submodule 
$F\cong \langle x^{\sigma}\rangle_F\subset A(n)$.
Moreover, it is also easy to see that the Lie algebra homomorphism 
$\sigma: H(n)_{\geq 0}\to F$ given by
$\sigma(x):={\rm tr}({\rm ad}_{H(n)/H(n)_{\geq 0}} x)$ for $x\in H(n)_{\geq 0}$
(see \cite[Pag. 155]{FS91}) is trivial.

 Moreover, using Lemma \ref{H-commutators}, 
it is straithforward to check that, in the notation of \cite{FS91}, we have  
the equality
$$[H(n)_{\geq 0}, H(n)_{\geq 0}]:=
H(n)_{\geq 0}^{(1)}=H(n)_{\geq 0}.$$
Therefore, using \cite[Thm. 3.6(2)]{FS91}, we get that 
$$H^2(H(n), A(n))=H^2(H(n)_{\geq 0}, F)\oplus \bigwedge^2 (H(n)/H(n)_{\geq 0}).$$
Since $\dim_F \bigwedge^2 (H(n)/H(n)_{\geq 0})=\binom{n}{2}$, we conclude using 
Proposition \ref{H(n)0-trivial} below.
\end{proof}

\begin{pro}\label{H(n)0-trivial}
The second cohomology group of $H(n)_{\geq 0}$ with coefficients in the trivial
module $F$ is given by
$$H^2(H(n)_{\geq 0},F)=
\begin{sis}
& \bigoplus_{i=1}^n \langle \overline{{\rm Sq}(x_i)} \rangle_F \bigoplus_{i=1}^n 
\langle \Omega_i\rangle \rangle_F \bigoplus \langle 
\Phi \rangle_F \text{ if } n\not\equiv -4\mod p,\\
&\bigoplus_{i=1}^n \langle \overline{{\rm Sq}(x_i)} \rangle_F \bigoplus_{i=1}^n 
\langle \Omega_i\rangle \rangle_F \bigoplus \langle 
\Phi \rangle_F\bigoplus \langle \Delta\rangle_F \text{ otherwise. }
\end{sis}$$
where $\Omega_i$ and $\Delta$ are the cocycles of Proposition \ref{H-triv-2}, 
$\Phi$ and $\Pi_{ij}$ (with 
$j\neq i'$) are the cocycles of Theorem \ref{H-finaltheorem} and 
$\overline{{\rm Sq}(x_i)}$ is the projection of ${\rm Sq}(x_i)$ onto 
$\langle 1\rangle_F\cong F $.
\end{pro}

\begin{proof}
We prove first that $$H^2(H(n)_{\geq 0}, F)=H^2(H(n)_{\geq 1}, F)^{H(n)_0},$$
where $H(n)_0$ acts on $H(n)_{\geq 1}$ via adjoint action. To this aim, 
consider the Hochschild-Serre spectral sequence with respect to the ideal
$H(n)_{\geq 1} \lhd H(n)_{\geq 0}$:
\begin{equation}\label{H-HS-A2}
E_2^{r,s}=H^r(H(n)_0,H^s(H(n)_{\geq 1},F))\Rightarrow H^{r+s}(H(n)_{\geq 0},F).
\end{equation}
By the Lemma \ref{H-commutators} and by homogeneity, 
it follows that
$$E_2^{1,1}=H^1(H(n)_0,H^1(H(n)_{\geq 1},F))=H^1(H(n)_0,C^1(H(n)_1,F))_{\0}=0.$$
Therefore we are left with showing that $H^2(H(n)_{\geq 0},F)=0$. 
First of all, we prove that $Z^1(H(n)_0,F)_{\0}=0$. 
Indeed, a homogeneous element 
$g\in C^1(H(n)_0,1)_{\0}$ can only take the following non-zero values 
$g(x_ix_{i'})=\alpha_i\cdot 1$, with $\alpha_i=\alpha_{i'}\in F$. 
The vanishing of $g$ follows from the following cocycle condition
\begin{equation*}
0=\d g(x_i^2,x_{i'}^2)=-4\sigma(i)g(x_ix_{i'})=-4\sigma(i)\alpha_i.\tag{*}
\end{equation*}
Consider now a homogeneous
cochain $f\in C^2(H(n)_0,F)_{\0}$. By applying the cocycle condition to the elements of $T_H$
and using homogeneity, one gets that $f_{|T_H}\in Z^1(H(n)_0,1)_{\0}$ 
which vanishes as proved above. Moreover, by adding to $f$ a coboundary, we can suppose that 
$f(x_i^2,x_{i'}^2)=0$
(see equation (*) of above). Therefore the only non-zero values of $f$
can be $f(x_ix_j,x_{i'}x_{j'})=\alpha_{ij}\cdot 1$ (for $j\neq i,i'$)  
with the obvious relations
$\alpha_{ij}=\alpha_{ji}$ and $\alpha_{ij}=-\alpha_{i'j'}$. We conclude by mean of the following 
cocycle condition
$$0=\d f(x_i^2,x_{i'}x_j,x_{i'}x_{j'})=-2\sigma(i)\alpha_{ij}+2\sigma(i)\alpha_{ij'},$$
which gives $\alpha_{ij}=\alpha_{ij'}=\alpha_{i'j'}=-\alpha_{ij}$ and hence $\alpha_{ij}=0.$

In order to compute $H^2(H(n)_{\geq 1},F)^{H(n)_0}$, we will use the same
strategy of Proposition \ref{K-H^2},  
that is to compute, step by step as $d$ increases, the truncated invariant cohomology groups
$$H^2\left(\frac{H(n)_{\geq 1}}{H(n)_{\geq d+1}},F\right)^{H(n)_0}.$$
By using the Hochschild-Serre spectral sequence associated to the ideal 
\begin{equation*}
H(n)_d=\frac{H(n)_{\geq d}}{H(n)_{\geq d+1}}\lhd \frac{H(n)_{\geq 1}}{H(n)_{\geq d+1}}\: ,
\end{equation*}
we obtain the same diagram as in \cite[Prop. 3.10]{Viv} (the vanishing of $E_2^{0,2}$ 
and the injectivity of the map $\alpha$ are proved in exactly the same way) and then
we take the cohomology with respect to $H(n)_0$.
An easy inspection of their proof shows that the Lemmas \ref{K-inv-(1,1)},
\ref{K-inv-(0,1)} and  \ref{K-H^1-(0,1)} of the preceding section (for the algebra $K(2m+1)$)
can be easily adapted to the present case simply by ignoring the variable $x_{2m+1}$ . 
In particular we get that (for $d\geq 1$) 
$$\begin{aligned}
& C^1(H(n)_1\times H(n)_d, F)^{H(n)_0}=
\begin{cases}
\langle \Phi_2\rangle_F & \text{ if } d=1,\\
\langle \Psi_2 \rangle_F & \text{ if } d=n(p-1)-5,\\ 
\end{cases}\\
&C^1(H(n)_d,F)^{H(n)_0}=0,\\
& H^1(H(n)_0,C^1(H(n)_d,F))=\begin{cases}
\oplus \langle \overline{{\rm Sq}(x_i)} \rangle_F & \text{ if } d=p-2,\\
\oplus_{i=1}^n \langle \omega_i \rangle_F & \text{ if } d=n(p-1)-p-2.\\
\end{cases}\\
\end{aligned}$$
where $\Phi_2$, $\Psi_2$ and $\omega_i$ are defined as in the case of $K(n)$ but ignoring 
the part involving the variable $x_{2m+1}=x_n$.

By definition $\overline{{\rm Sq}(x_i)}$ is the restriction of ${\rm Sq}(x_i)$ and it is easy to see
that $\omega_i$ is the restriction of $\Omega_i$. Moreover if we extend 
$\Phi_2$ by $0$ outside $H(n)_1\times H(n)_1$, then it is clear that $\Phi_2\in 
H^2(H(n)_{\geq 0},F)\subset H^2(H(n)_{\geq 1},F)^{H(n)_0}$ and that $\Phi_2$
is the restriction of the cocycle $\Phi$ (see also \cite[Prop. 3.7]{Viv}). 

Finally, suppose that there is a lifting of $\Psi_2$ to a global $H(n)_0$-invariant 
cocycle of $Z^2(H(n)_{\geq 1},$ $A(n))$, which we will continue to call $\Psi_2$. 
Then using the cocycle condition $0=\d \Psi_{2_{|H(n)_1}}$
together with Lemma \ref{H-commutators} and proceeding by induction on the degree, 
it is easy to see that $\Psi_2$ must agree with $\Delta$ on the couples 
$(x^a,x^b)\in H(n)_{\geq 1}\times H(n)_{\geq 1}$ such that $a+b=\sigma$ and we 
know from Proposition \ref{H-triv-2} that $\Delta$ is an antisymmetric cocycle 
if and only if $n\equiv -4 \mod p$.

\end{proof}

\begin{lem}\label{H-commutators}
Let $d$ be an integer greater or equal to $-1$. 
Then $$[H(n)_1,H(n)_d]=H(n)_{d+1}.$$
\end{lem}
\begin{proof}
The proof is the same as the first part of Lemma \ref{K-commutators} (where we
consider elements belonging to $A(2m)\subset K(2m+1)$) except for the fact that 
we do not have to consider the elements $x_i$ because they have degree $-1$ 
and the element $x^{\sigma}$ which does not belong to $H(n)$.
\end{proof}


\begin{thebibliography}{TESI2}

\addcontentsline{toc}{section}{References}

\bibitem[BW88]{BW} R.E. Block and R.L. Wilson: \emph{Classification of the restricted simple
Lie algebras}.  J. Algebra  114  (1988), 115--259.



\bibitem[CE48]{CE} C. Chevalley and S. Eilenberg: \emph{Cohomology theory of Lie groups
and Lie algebras}. Trans. Amer. Math. Soc. 63, (1948), 85--124.

\bibitem[DG70]{DG} M. Demazure and P. Gabriel: \emph{Groupes alg\'ebriques}.
Tome I: G\'eom\'etrie alg\'ebrique, g\'en\'eralit\'es, groupes commutatifs (French).
Avec un appendice  Corps de classes local par Michiel Hazewinkel.
Masson and Cie, Editeur, Paris; North-Holland Publishing Co., Amsterdam, 1970.


\bibitem[FAR86]{Far} R. Farnsteiner: \emph{Central extensions and invariant forms of 
graded Lie algebras}. Algebras, Groups, Geom. 3 (1986), 431--455.

\bibitem[FS88]{FS} R. Farnsteiner and H. Strade: \emph{Modular Lie algebras
and their representation}. Monographs and textbooks in pure and applied mathematics,
vol. 116 (Dekker, New York, 1988).

\bibitem[FS91]{FS91} R. Farnsteiner and H. Strade: \emph{Shapiro's Lemma and its consequences 
in the cohomology theory of modular Lie algebras}. Math. Zeit. 206 (1991), 153--168.

\bibitem[GER64]{GER1} M. Gerstenhaber: \emph{On the deformation of rings and algebras}.  
Ann. of Math. 79  (1964), 59--103.


\bibitem[HS53]{HS} G. Hochschild and J.-P. Serre: \emph{Cohomology of Lie algebras}.
Ann. of Math. 57 (1953), 591--603.

\bibitem[KS66]{KS} A. I. Kostrikin and I. R. Shafarevich: \emph{Cartan's pseudogroups
and the $p$-algebras of Lie} (Russian).  Dokl. Akad. Nauk SSSR  168  (1966), 740--742.
English translation: Soviet Math. Dokl. 7 (1966), 715--718.


\bibitem[MEL80]{MEL} G. M. Melikian: \emph{Simple Lie algebras of characteristic $5$}
(Russian).  Uspekhi Mat. Nauk  35  (1980), no. 1 (211), 203--204.




\bibitem[PS01]{PS3} A. Premet and  H. Strade: \emph{Simple Lie algebras of small characteristic.
III. The toral rank 2 case}.  J. Algebra  242  (2001), 236--337.


\bibitem[RUD71]{RUD} A. N. Rudakov: \emph{Deformations of simple Lie algebras} (Russian).
 Izv. Akad. Nauk SSSR Ser. Mat. 35 (1971), 1113--1119.

\bibitem[SEL67]{SEL} G. B. Seligman: \emph{Modular Lie algebras}. Ergebnisse der Mathematik 
und ihrer Grenzgebiete, Band 40 Springer-Verlag, New York 1967. 






\bibitem[STR98]{STR6} H. Strade: \emph{The classification of the simple modular Lie algebras.
VI. Solving the final case}.  Trans. Amer. Math. Soc.  350  (1998), 2553--2628.


\bibitem[STR04]{STR} H. Strade: \emph{Simple Lie algebras over fields of 
positive characteristic I:
Structure theory}. De Gruyter Expositions in Mathematics, 38 (Walter de Gruyter, Berlin, 2004).

\bibitem[VIV1]{Viv} F. Viviani: \emph{Infinitesimal deformations of restricted 
simple Lie algebras I}. To appear on J. Algebra 
(available at arXiv:math.RA/0612861).

\bibitem[VIV2]{Viv2} F. Viviani: \emph{Deformations of the restricted Melikian algebra}.
To appear on Comm. Algebra (available at arXiv:math.RA/0702594).

\bibitem[VIV3]{Viv3} F. Viviani: \emph{Deformations of simple finite group schemes}, submitted.
(available at arXiv:0705.0821).






\end{thebibliography}
\end{document}